\documentclass[]{article}
\usepackage{amsmath}
\usepackage{amsthm}
\usepackage{url}
\usepackage{authblk}

\usepackage[dvips]{graphicx}
\usepackage{color}
\usepackage{amsfonts,amssymb,amsmath,}
\usepackage[T1]{fontenc}
\usepackage{ae}

\newtheorem{thm}{Theorem}

\newtheorem{coro}{Corollary}[thm] 

\newtheorem{lem}{Lemma}        

\numberwithin{equation}{section} 

\allowdisplaybreaks[4]

\begin{document}
\title{The distribution of $k$-free numbers and the derivative of the Riemann zeta-function}
\author[]{Xianchang Meng}
\date{}
\maketitle
\begin{abstract}
Under the Riemann Hypothesis, we connect the distribution of $k$-free numbers with the derivative of the Riemann zeta-function at nontrivial zeros of $\zeta(s)$. Moreover, with additional assumptions, we prove the existence of a limiting distribution of $e^{-\frac{y}{2k}}M_k(e^y)$ and study the tail of the limiting distribution, where $M_k(x)=\sum_{n\leq x}\mu_k(n)-\frac{x}{\zeta(k)}$ and $\mu_k(n)$ is the characteristic function of $k$-free numbers. Finally, we make a conjecture about the maximum order of $M_k(x)$ by heuristic analysis on the tail of the limiting distribution.
\end{abstract}

\let\thefootnote\relax
\let\thefootnote\relax

\setcounter{tocdepth}{1}

\section{Introduction and Results}

Let $\mu_k(n)$ be the characteristic function of $k$-free numbers, where $k\geq 2$.
 Let $\widetilde{M}_k(x)=\sum_{n\leq x}\mu_k(n)$ be the number of $k$-free integers $\leq x$, and $M_k(x)=\widetilde{M}_k(x)-\frac{x}{\zeta(k)}$.
Using elementary arguments, one can derive
$$\widetilde{M}_k(x)=\frac{x}{\zeta(k)}+O(x^{\frac{1}{k}}).$$
Many authors have worked to improve the error term. The best unconditional result is due to Walfisz \cite{Walf},
$$M_k(x)\ll x^{\frac{1}{k}}\exp\{-ck^{-\frac{8}{5}}\log^{\frac{3}{5}}x\log\log^{-\frac{1}{5}}x\},$$
where $c>0$ is an absolute constant.
In the opposite direction, Evelyn and Linfoot \cite{Evelyn} proved that
$$M_k(x)=\Omega(x^{\frac{1}{2k}}).$$

Under the Riemann Hypothesis, Montgomery and Vaughan \cite{MontVaug2} showed that
$$M_k(x)\ll x^{\frac{1}{k+1}+\epsilon}, \quad \forall \epsilon >0.$$
Later, several authors made contributions to the improvement of the error term under the Riemann Hypothesis, such as, Graham \cite{Grah}, Baker and Pintz \cite{Baker-Pintz}, Jia \cite{Jia}, \cite{Jia2}, Graham and Pintz \cite{Grah-Pintz}, and Baker and Powell \cite{Baker-Powell}, etc. They improved Montgomery and Vaughan's result for various $k$. Their bounds have the shape
$M_k(x)\ll x^{E(k)}$, where $E(k)\sim \frac{1}{k}$ as $k\rightarrow\infty$. For $k=2$, the best result under the Riemann Hypthosis is due to Jia \cite{Jia2} in which the exponent is $\frac{17}{54}+\epsilon$.
For a nice survey, see \cite{Pappa} or \cite{Sand}. However, there is still a large gap to the conjectured result,
$$M_k(x)\ll x^{\frac{1}{2k}+\epsilon}, \quad \forall \epsilon >0.$$

In this paper, we take a different approach to studying the distribution of $k$-free numbers by connecting it to the analytic properties of $\zeta(s)$.
Define
$$J_{-l}(T)=\sum_{0<\gamma\leq T}\frac{1}{|\zeta'(\rho)|^{2l}},$$
where $l\in \mathbb{R}$, $\rho$ is a zero of $\zeta(s)$, and $\gamma$ is the imaginary part of $\rho$. For the existence of the above sum, we implicitly assume the zeros of $\zeta(s)$ are simple. Gonek \cite{Gonek} and Hejhal \cite{Hej} independently conjectured that
\begin{equation}\label{G-conj}
  J_{-l}(T)\asymp T(\log T)^{(l-1)^2}.
\end{equation}
For $l=1$, Gonek \cite{Gonek} proved that $J_{-1}(T)\gg T$ subject to the Riemann Hypothesis and the simplicity of zeros. Moreover, he conjectured in \cite{Gonek2} that
\begin{equation}\nonumber
  J_{-1}(T)\sim \frac{3}{\pi^2}T.
\end{equation}
Hughes, Keating, and O'Connell \cite{HKC} used a different heuristic method based on random matrix theory and made the conjecture
\begin{equation}\nonumber
  J_{-l}(T)\sim C_l \frac{T}{2\pi}\left(\log \frac{T}{2\pi}\right)^{(l-1)^2} \quad\mbox{for}\quad l<\frac{3}{2},
\end{equation}
where $C_l$ is a constant depending on $l$.

Recently, Ng \cite{Ng} connected the summatory function of the M\"{o}bius function, $M(x)=\sum_{n\leq x}\mu(n)$, with the behavior of $J_{-1}(T)$. He showed that, under the Riemann Hypothesis, $J_{-1}(T)\ll T$ implies the so called weak Mertens conjecture, which asserts
\begin{equation}\label{weak-M-conj}
\int_2^X\left(\frac{M(x)}{x}\right)^2dx\ll \log X.
\end{equation}
He also proved that $e^{-y/2}M(e^y)$ has a limiting distribution under the assumptions of the Riemann Hypothesis and $J_{-1}(T)\ll T$, and studied the tail of the distribution using Montgomery's probabilistic methods. On the other hand, under the Riemann Hypothesis, Titchmarsh \cite{Titch} (Chapter XIV, page 376-380) showed that the weak Mertens conjecture implies the simplicity of zeros of $\zeta(s)$, the estimate $\frac{1}{\zeta'(\rho)}=O(|\rho|)$, and the convergence of the series $\sum_{\rho}\frac{1}{|\rho\zeta'(\rho)|^2}$.

Motivated by Ng's work and Titchmarsh's argument, we connect the estimation of $J_{-1}(T)$ to the distribution of $k$-free numbers. In our work, we don't require as strong an assumption $J_{-1}(T)\ll T$ as used in \cite{Ng}.
 Under the Riemann Hypothesis, we show that
\begin{equation}\label{GH}
  J_{-1}(T)=\sum_{0<\gamma\leq T}\frac{1}{|\zeta'(\rho)|^2}\ll_{\epsilon} T^{1+\epsilon} \quad \forall \epsilon>0,
\end{equation}
and
 \begin{equation}\label{equiv-M}
    \int_1^X \left(\frac{M_k(x)}{x^{\frac{1}{2k}}}\right)^2\frac{dx}{x}\ll_k \log X \quad\forall k\geq 2,
  \end{equation}
 are equivalent. We can view (\ref{equiv-M}) as an analogue of the weak Mertens conjecture. Further, under the Riemann Hypothesis and (\ref{GH}), we prove the following results about the distribution of $k$-free numbers.

\begin{thm}\label{thm-average}
  Under the Riemann Hypothesis, (\ref{GH}) and (\ref{equiv-M}) are equivalent.
\end{thm}
\noindent\textbf{Remark 1.} One can use (\ref{GH}), and Lemmas \ref{lemma-M} and \ref{lem-square}, to show that $M_k(x)\ll_{k, \epsilon} x^{\frac{1}{2k}}(\log x)^{\frac{1}{2}-\frac{1}{2k}+\epsilon}$
except on a set of finite logarithmic measure.

\begin{thm}\label{thm-asymp}
  Assume the Riemann Hypothesis and (\ref{GH}). Then, for any $k\geq 2$, we have
  $$\int_1^{X}\left(\frac{M_k(x)}{x^{\frac{1}{2k}}}\right)^2\frac{dx}{x}\sim \beta_k \log X,$$
  where
  $$\beta_k=\sum_{\gamma>0}\frac{2|\zeta(\frac{\rho}{k})|^2}{|\rho\zeta'(\rho)|^2}.$$
\end{thm}

The main part of our proof is to show that $\phi(y)=e^{-\frac{y}{2k}}M_k(e^{y})$ is a $B^2$-\textit{almost periodic} function, which means, for any $\epsilon>0$, there exists a real valued trigonometric polynomial $P_{N(\epsilon)}(y)=\sum_{n=1}^{N(\epsilon)} r_n(\epsilon)e^{i\lambda_n(\epsilon)y}$, such that
\begin{equation}\label{B2-almost}
\limsup_{Y\rightarrow\infty}\frac{1}{Y}\int_0^Y |\phi(y)-P_{N(\epsilon)}(y)|^2 dy<\epsilon^2.
\end{equation}
By the work of Besicovitch (see \cite{Besi}, Chapter II of \cite{Besi2}, or Theorem 1.14 of \cite{A-N-S}), a Parseval type identity is true for $B^2$-almost periodic functions. Moreover, $B^2$-almost periodic functions possess limiting distributions (see Theorem 2.9 of \cite{A-N-S}). Thus, we get the following result.

\begin{thm}\label{thm-distri}
  Assume the Riemann Hypothesis and (\ref{GH}). Then, for any $k\geq 2$,  $e^{-\frac{y}{2k}}M_k(e^{y})$ has a limiting distribution $\nu=\nu_k$ on $\mathbb{R}$, that is
  $$\lim_{Y\rightarrow\infty}\frac{1}{Y}\int_{0}^{Y} f(e^{-\frac{y}{2k}}M_k(e^{y}))dy=\int_{-\infty}^{\infty}f(x)d\nu(x)$$
  for all bounded Lipschitz continuous functions $f$ on $\mathbb{R}$.
\end{thm}

The Linear Independence Conjecture (LI) states that, under the Riemann Hypothesis, the positive imaginary ordinates of the zeros of $\zeta(s)$ are linearly independent over $\mathbb{Q}$. If we add the assumption of LI, we can show the following corollaries. The first corollary is similar to Theorem 1.9 of \cite{A-N-S}. We omit the proof.
\begin{coro}\label{cor-fouri}
  Assume the Riemann Hypothesis, (\ref{GH}), and LI. Then the Fourier transform $\hat{\nu}(\xi)=\int_{-\infty}^{\infty}e^{-i\xi t}d\nu(t)$ exists and equals
  $$\widehat{\nu}(\xi)=\prod_{\gamma>0}\widetilde{J}_0\left(\frac{2|\zeta(\frac{\rho}{k})|\xi}{|\rho\zeta'(\rho)|}\right),$$
  where $\widetilde{J}_0(z)$ is the Bessel function
  $\widetilde{J}_0(z)=\sum_{m=0}^{\infty}\frac{(-1)^m(\frac{1}{2}z)^{2m}}{(m!)^2}$.
\end{coro}

\begin{coro}\label{app-LI}
  Assume the Riemann Hypothesis, (\ref{GH}), and LI.
Then, for any $k\geq 2$ and any $\epsilon>0$, we have
\begin{equation*}
\exp\left(-\widetilde{c}_1 V^{\frac{2k}{k-1}+\epsilon}\right)\leq \nu([V, \infty))\leq \exp\left(-\widetilde{c}_2 V^{\frac{2k}{k-1}-\epsilon}\right),
\end{equation*}
for some constants $\widetilde{c}_1, \widetilde{c}_2 >0$ depending on $k$ and $\epsilon$.
\end{coro}

Under the Riemann Hypothesis, (\ref{G-conj}) (for $l< \frac{3}{2}$), and LI, with a refined analysis, we can prove a more precise upper bound for the tail of the limiting distribution. Further analysis of the bounds for $\nu([V, \infty))$ suggests the following conjecture.\\
\textbf{Large deviation conjecture. } There exist positive constants $c'_1$, $c'_2$ such that for large $V$,
\begin{equation}\label{conj-order}
\exp\left(-c_2' \frac{V^{\frac{2k}{k-1}}}{(\log V)^{\frac{1}{2(k-1)}}}\right)\ll \nu([V, \infty))\ll \exp\left(-c_1' \frac{V^{\frac{2k}{k-1}}}{(\log V)^{\frac{1}{2(k-1)}}}\right).
\end{equation}

\begin{thm}\label{thm-la-dev-upper}
	Assume the Riemann Hypotheis, (\ref{G-conj}) (for $l< \frac{3}{2}$), and LI. Then, there exists a constant $c_1''>0$ such that
	\begin{equation}\label{conj-up-bound}
	\nu([V, \infty))\ll \exp\left(-c_1'' \frac{V^{\frac{2k}{k-1}}}{(\log V)^{\frac{1}{2(k-1)}+o(1)}}\right).
	\end{equation}
\end{thm}

\noindent\textbf{Remark 2.} Bounds for the tail of the probabilistic measure can be used to heuristically estimate the maximum variation of $M_k(x)$.

With the above Large deviation conjecture and a similar heuristic argument to section 4.3 of \cite{Ng}, we make the following conjecture.

\noindent\textbf{Conjecture.}  For any $k\geq 2$, there exists a number $C=C_k>0$, such that
\begin{equation}\label{M-conj}
\overline{\underline{\lim}}_{x\rightarrow\infty}\frac{M_k(x)}{x^{\frac{1}{2k}}(\log\log x)^{\frac{k-1}{2k}}(\log\log\log x)^{\frac{1}{4k}}}=\pm C_k.
\end{equation}

For the proof of Theorem \ref{thm-la-dev-upper}, we need a result on the moments of the Riemann zeta-function.
\begin{thm}\label{thm-moments-zeta}
	Assume the Riemann Hypothesis. For any fixed integer $l\geq 1$ and $0<w<1$,
	\begin{equation}\label{thm-moms-one}
	\sum_{0<\gamma\leq T}|\zeta(1-w\rho)|^{2l}= C_{w,l}T\log T+O_{w, l}(T(\log T)^{\frac{1}{2}}),
	\end{equation}
	and
	\begin{equation}\label{thm-moms-two}
	\sum_{0<\gamma\leq T}\frac{1}{|\zeta(1-w\rho)|^{2l}}= C'_{w,l}T\log T+O_{w, l}(T(\log T)^{\frac{1}{2}}),
	\end{equation}
	where
	$$C_{w,l}=\frac{1}{2\pi}\sum_{n=1}^{\infty}\frac{d_l^2(n)}{n^{2-w}}, \   C'_{w,l}=\frac{1}{2\pi}\sum_{n=1}^{\infty}\frac{{\tilde{d}}_l^2(n)}{n^{2-w}},$$
    $d_l(n)$ denotes the number of ways $n$ may be written as a product of $l$ factors, and ${\tilde{d}}_l(n)=(\underbrace{\mu*\mu*\cdots*\mu}_{l~\text{times}})(n)$.
\end{thm}

\noindent\textbf{Remark 3.} Using our techniques, we may replace the assumption of $J_{-1}(T)\ll T$ in \cite{Ng} by $J_{-1}(T)\ll T^{2-\epsilon}$ for any $\epsilon>0$, and still obtain the conclusions in Theorem 1, parts ii), iii) and iv), Theorem 2, and Theorem 3 of \cite{Ng}.
Let $k=1$ and $\theta=2-\epsilon$ in our Lemma \ref{lem-square}, we get
 $$\int_{\log Z}^{\log Z+1}\left|\sum_{T\leq \gamma\leq X}\frac{1}{\rho\zeta'(\rho)}e^{i\gamma y}\right|^2 dy\ll \frac{1}{T^{\epsilon/2}}.$$
With this bound and a similar argument to Theorems \ref{thm-average}, \ref{thm-asymp}, and \ref{thm-distri} we have
\begin{thm}
	The Riemann Hypothesis and $J_{-1}(T)\ll T^{2-\epsilon}$ for any fixed $\epsilon>0$ imply
	\begin{description}
		\item[\quad 1)] the weak Mertens conjecture (\ref{weak-M-conj});
		\item[\quad 2)] $M(x)=\sum_{n\leq x}\mu(n)\ll x^{1/2}(\log\log x)^{3/2}$ except on a set of finite logarithmic measure;
		\item[\quad 3)] $$\int_1^{X}\left(\frac{M(x)}{x^{1/2}}\right)^2\frac{dx}{x}\sim \beta \log X,$$
		where
		$\beta=\sum_{\gamma>0}\frac{2}{|\rho\zeta'(\rho)|^2};$
		\item[\quad 4)] $e^{-y/2}M(e^y)$ has a limiting distribution on $\mathbb{R}$.
		
	\end{description}
\end{thm}

\section{Main Lemmas and Proofs}\label{sec-lemma}
Since (\ref{GH}) implies that the zeros of $\zeta(s)$ are simple, in this section, we implicitly assume the simplicity of zeros of $\zeta(s)$.
And the implicit constants in our estimates may depend only on $k$ or $\epsilon$, unless otherwise specified.

\begin{lem}\label{lem-G}
  Assume $J_{-1}(T)\ll T^{\theta}$ for some $\theta\geq 1$. Then, for $a>\frac{\theta+1}{2}, b>\theta$ and any $\epsilon>0$, we have
  $$\sum_{\gamma>T}\frac{1}{\gamma^a|\zeta'(\rho)|}\ll \frac{1}{T^{a-\frac{\theta+1}{2}-\epsilon}}, \quad \mbox{and}\quad\sum_{\gamma>T}\frac{1}{\gamma^b|\zeta'(\rho)|^2}\ll \frac{1}{T^{b-\theta}} .$$
\end{lem}

\noindent \emph{Proof.}
Let $N(T)$ be the number of zeros of $\zeta(s)$ in the region $0<\sigma<1$, $0<t\leq T$. By Theorem 1.7 in \cite{Ivic},
\begin{equation}\label{zero-count}
  N(T)=\frac{T}{2\pi}\log\frac{T}{2\pi}-\frac{T}{2\pi}+O(\log T).
\end{equation}
By the Cauchy-Schwarz inequality, and (\ref{zero-count}),
\begin{equation*}
J_{-\frac{1}{2}}(T)=\sum_{0<\gamma\leq T}\frac{1}{|\zeta'(\rho)|}\ll \left(\sum_{0<\gamma\leq T}\frac{1}{|\zeta'(\rho)|^2}\right)^{\frac{1}{2}}\left(\sum_{0<\gamma\leq T} 1\right)^{\frac{1}{2}}\ll \frac{(\log T)^{\frac{1}{2}}}{T^{a-\frac{\theta+1}{2}}}
\end{equation*}
Let $f(t)=\frac{1}{t^a}$. Then, $f'(t)=-at^{-a-1}$, and by partial summation,
$$\sum_{\gamma>T}\frac{1}{\gamma^a|\zeta'(\rho)|}=f(t)J_{-\frac{1}{2}}(t)|_{T}^{\infty}-\int_T^{\infty}J_{-\frac{1}{2}}(t)f'(t)dt\ll \frac{1}{T^{a-\frac{\theta+1}{2}-\epsilon}}.$$
Similarly, using $J_{-1}(T)\ll T^\theta$ and partial summation, we can prove the second formula. \qed

\vspace{1em}
We also need the following lemma from \cite{Ng}.
\begin{lem}[\cite{Ng}, Lemma 3]\label{lemma3-Ng}
  Assume the Riemann Hypothesis. There exists a sequence of numbers $\mathcal{T}=\{T_n\}_{n=0}^{\infty}$ which satisfies
  $$n\leq T_n\leq n+1 \quad \mbox{and}\quad \frac{1}{\zeta(\sigma+iT_n)}=O(T_n^{\epsilon})\quad (-1\leq\sigma\leq 2).$$
\end{lem}

One can use classical tools, like Perron's formula, contour integration, and functional equation, to prove the following lemma. As the argument is standard, we do not provide details of the proof.
\begin{lem}\label{lemma-3}
  Assume the Riemann Hypothesis and that all zeros of $\zeta(s)$ are simple. For $T\in \mathcal{T}$,
  $$\widetilde{M}_k(x)=\sum_{n\leq x}\mu_k(n)=\frac{x}{\zeta(k)}+\sum_{|\gamma|<T}\frac{\zeta(\frac{\rho}{k})}{\rho\zeta'(\rho)}x^{\frac{\rho}{k}}+\widetilde{E}(x, T), $$
  where
  $$\widetilde{E}(x, T)\ll_{k, \epsilon} \frac{x\log x}{T}+\frac{x}{T^{\frac{1}{2}-\epsilon}}+1.$$
\end{lem}

The two Lemmas below are the main lemmas we use to prove Theorems \ref{thm-average}, \ref{thm-asymp}, and \ref{thm-distri}. We will give their proofs in the following subsections.

\begin{lem}\label{lemma-M}
  Assume the Riemann Hypothesis and (\ref{GH}). For $x\geq 2$, $T\geq 2$, and any $\epsilon>0$,
  $$\widetilde{M}_k(x)=\frac{x}{\zeta(k)}+\sum_{|\gamma|<T}\frac{\zeta(\frac{\rho}{k})}{\rho\zeta'(\rho)}x^{\frac{\rho}{k}}+E(x, T),$$
  where
  $$E(x, T)\ll_{k,\epsilon} \frac{x\log x}{T}+\frac{x}{T^{\frac{1}{2}-\epsilon}}+\frac{x^{\frac{1}{2k}}}{T^{\frac{1}{2k}-\epsilon}}+1.$$
\end{lem}

\begin{lem}\label{lem-square}
 Let $k\geq 1$ be an integer. Assume the Riemann Hypothesis and $J_{-1}(T)\ll T^{\theta}$ for some $1\leq\theta<1+\frac{1}{k}$. Then for any $\epsilon>0$,  $Z\geq 0$ and $0<T<X$,
  $$\int_{\log Z}^{\log Z+1}\left|\sum_{T\leq \gamma\leq X}\frac{u_{k}(\rho)}{\rho\zeta'(\rho)}e^{\frac{i\gamma y}{k}}\right|^2 dy\ll_{k, \epsilon} \frac{1}{T^{1+\frac{1}{k}-\theta-\epsilon}},$$
  where $u_k(\rho)=\zeta(\frac{\rho}{k})$ if $k\geq 2$, $u_k(\rho)=1$ if $k=1$.
\end{lem}

\subsection{Proof of Lemma \ref{lemma-M}}
Let $T\geq 2$ and $n\leq T\leq n+1$. Without loss of generality, assume $n\leq T_n\leq T\leq n+1$. Then, by Lemma \ref{lemma-3},
\begin{equation}\nonumber
 \widetilde{M}_k(x)=\frac{x}{\zeta(k)}+\sum_{|\gamma|\leq T}\frac{\zeta(\frac{\rho}{k})}{\rho\zeta'(\rho)}x^{\frac{\rho}{k}}-\sum_{T_n\leq|\gamma|\leq T}\frac{\zeta(\frac{\rho}{k})}{\rho\zeta'(\rho)}x^{\frac{\rho}{k}}+\widetilde{E}(x, T).
\end{equation}
By the Cauchy-Schwarz inequality,
\begin{equation}\nonumber
\left|\sum_{T_n\leq\gamma\leq T}\frac{\zeta(\frac{\rho}{k})}{\rho\zeta'(\rho)}x^{\frac{\rho}{k}}\right| \leq x^{\frac{1}{2k}}\left(\sum_{T_n\leq \gamma\leq T}\frac{1}{|\zeta'(\rho)|^2}\right)^\frac{1}{2}\cdot\left(\sum_{T_n\leq \gamma\leq T}\left|\frac{\zeta(\frac{\rho}{k})}{\rho}\right|^2\right)^{\frac{1}{2}}.
\end{equation}
Under the Riemann Hypothesis, by section 13.1 in \cite{Titch},
\begin{equation}\label{zeta-est}
  |\zeta(\sigma+it)|\ll(|t|+2)^{\frac{1}{2}-\sigma+\epsilon}, ~\mbox{for}~0\leq \sigma\leq \frac{1}{2}.
\end{equation}
By (\ref{zeta-est}), we have
\begin{equation}\label{zeta-rho}
  \left|\zeta(\frac{\rho}{k})\right|\ll_{k, \epsilon} |\gamma|^{\frac{1}{2}-\frac{1}{2k}+\epsilon},
\end{equation}
Thus, by (\ref{GH}), (\ref{zeta-rho}), and (\ref{zero-count}),
\begin{equation}\nonumber
\left|\sum_{T_n\leq\gamma\leq T}\frac{\zeta(\frac{\rho}{k})}{\rho\zeta'(\rho)}x^{\frac{\rho}{k}}\right|\ll_{k, \epsilon} x^{\frac{1}{2k}}T^{\frac{1}{2}+\epsilon}\frac{1}{T^{\frac{1}{2}+\frac{1}{2k}-\epsilon}}\ll_{k, \epsilon} \frac{x^{\frac{1}{2k}}}{T^{\frac{1}{2k}-\epsilon}}.
\end{equation}
The desired result follows.\qed

\subsection{Proof of Lemma \ref{lem-square}}
In this lemma, we assume the bound
\begin{equation}\label{lem-M-assump}
J_{-1}(T)\ll T^{\theta},
\end{equation}
for some $1\leq \theta<1+\frac{1}{k}$, and this implies, by the Cauchy-Schwarz inequality, that
\begin{equation}\label{J-1/2-est-pf-lem1}
J_{-\frac{1}{2}}(T)=\sum_{0<\gamma\leq T}\frac{1}{|\zeta'(\rho)|}\ll T^{\frac{\theta+1}{2}} (\log T)^{\frac{1}{2}}.
\end{equation}
We have
\begin{eqnarray}
 && \int_{\log Z}^{\log Z+1}\left|\sum_{T\leq \gamma\leq X}\frac{u_k(\rho)}{\rho\zeta'(\rho)}e^{\frac{i\gamma y}{k}}\right|^2 dy
  =\sum_{T\leq \gamma\leq X}\sum_{T\leq \gamma'\leq X}\frac{u_k(\rho)\overline{u_k(\rho')}}{\rho\zeta'(\rho)\overline{\rho'\zeta'(\rho')}}\int_{\log Z}^{\log Z+1} e^{\frac{i(\gamma-\gamma')y}{k}}dy\nonumber\\
  &&\quad\quad\ll_k \sum_{T\leq \gamma\leq X}\sum_{T\leq \gamma'\leq X}\left|\frac{u_k(\rho)\overline{u_k(\rho')}}{\rho\zeta'(\rho)\overline{\rho'\zeta'(\rho')}}\right|\min\left(1, \frac{1}{|\gamma-\gamma'|}\right).\nonumber
\end{eqnarray}
Note that $\rho=\frac{1}{2}+i\gamma$ and $\rho'=\frac{1}{2}+i\gamma'$ denote zeros of $\zeta(s)$. We break the last sum into two parts:
$$\Sigma_1: \mbox{~the sum over~} |\gamma-\gamma'|\leq 1,
\mbox{\quad and\quad}
\Sigma_2: \mbox{~the sum over~} |\gamma-\gamma'|>1.$$

In the following proof, we always assume $T\leq \gamma'\leq X$. Without loss of generality, we assume  $0<\epsilon<\frac{1+\frac{1}{k}-\theta}{10}$. By (\ref{zeta-rho}), for any $\epsilon>0$,
\begin{equation}\label{pf-lem-M-U-est}
|u_k(\rho)|\ll |\gamma|^{\frac{1}{2}-\frac{1}{2k}+\epsilon}.
\end{equation}
For the first sum, noting that $1-\frac{1}{k}-2\epsilon>\theta$, by the Cauchy-Schwarz inequality, (\ref{lem-M-assump}), and Lemma \ref{lem-G},
\begin{eqnarray}
  \Sigma_1&\ll&\sum_{T\leq \gamma\leq X}\frac{|u_k(\rho)|}{|\rho\zeta'(\rho)|}\sum_{\gamma-1\leq\gamma'\leq\gamma+1}\frac{|u_k(\rho')|}{|\rho'\zeta'(\rho')|}\nonumber\\
  &\ll&\left(\sum_{T\leq \gamma\leq X}\frac{|u_k(\rho)|^2}{|\rho\zeta'(\rho)|^2}\right)^{\frac{1}{2}}\left(\sum_{T\leq \gamma\leq X}\left(\sum_{\gamma-1\leq\gamma'\leq\gamma+1}\frac{|u_k(\rho')|}{|\rho'\zeta'(\rho')|}\right)^2\right)^{\frac{1}{2}}\nonumber\\
  &\ll&\left(\sum_{T\leq \gamma\leq X}\frac{|u_k(\rho)|^2}{|\rho\zeta'(\rho)|^2}\right)^{\frac{1}{2}}\left(\sum_{T\leq \gamma\leq X}\left(\sum_{\gamma-1\leq\gamma'\leq\gamma+1}\frac{|u_k(\rho')|^2}{|\rho'\zeta'(\rho')|^2}\right)\cdot\log\gamma\right)^{\frac{1}{2}}\nonumber\\
  &\ll&\left(\sum_{T\leq \gamma\leq X}\frac{1}{|\gamma|^{1+\frac{1}{k}-2\epsilon}|\zeta'(\rho)|^2}\right)^{\frac{1}{2}}\left(\sum_{T\leq \gamma\leq X}\sum_{\gamma-1\leq\gamma'\leq\gamma+1}\frac{|u_k(\rho')|^2}{|\rho'\zeta'(\rho')|^2}\log\gamma'\right)^{\frac{1}{2}}\nonumber\\
  &\ll&\frac{1}{T^{\frac{1}{2}+\frac{1}{2k}-\frac{\theta}{2}-\epsilon}}\left(\sum_{T\leq \gamma'\leq X}m(\gamma')\frac{|u_k(\rho')|^2}{|\rho'\zeta'(\rho')|^2}\cdot\log\gamma'\right)^{\frac{1}{2}},\nonumber
\end{eqnarray}
where $m(\gamma')=\#\{\gamma: \gamma-1\leq \gamma'\leq \gamma+1\}\ll \log \gamma'$ by (\ref{zero-count}). Since $1+\frac{1}{k}-3\epsilon>\theta$, by Lemma \ref{lem-G} and (\ref{pf-lem-M-U-est}),
\begin{equation}\label{sum1}
  \Sigma_1\ll \frac{1}{T^{\frac{1}{2}+\frac{1}{2k}-\frac{\theta}{2}-\epsilon}}\left(\sum_{T\leq \gamma'\leq X}\frac{1}{|\gamma'|^{1+\frac{1}{k}-3\epsilon}|\zeta'(\rho')|^2}\right)^{\frac{1}{2}} \ll\frac{1}{T^{1+\frac{1}{k}-\theta-3\epsilon}}.
\end{equation}

We write $\Sigma_2$ as follows,
\begin{equation}\label{sum2}
  \Sigma_2=\sum_{T\leq \gamma\leq X}\frac{|u_k(\rho)|}{|\rho\zeta'(\rho)|}\sum_{\substack{T\leq\gamma'\leq X\\ |\gamma-\gamma'|>1}}\frac{|u_k(\rho')|}{|\rho'\zeta'(\rho')||\gamma-\gamma'|}.
\end{equation}
Taking $N=[\frac{1}{\epsilon}]+2$, and let
$$1>a_1>a_2>a_3>\cdots>a_N=0.$$
Then, by (\ref{sum2}) we have
\begin{equation}\label{sum2-sig}
  \Sigma_2=\sum_{l=1}^{2N+1} \sigma_l,
\end{equation}
where $$\sigma_l=\sum_{T\leq \gamma\leq X}\frac{|u_k(\rho)|}{|\rho\zeta'(\rho)|}\sum_{\gamma'\in L_l}\frac{|u_k(\rho')|}{|\rho'\zeta'(\rho')||\gamma-\gamma'|},$$
and
$L_1:T\leq \gamma'<\gamma-\gamma^{a_1},$
$L_2:\gamma-\gamma^{a_1}\leq \gamma'<\gamma-\gamma^{a_2},$
$\cdots,$
$L_{N-1}:\gamma-\gamma^{a_{N-2}}\leq \gamma'<\gamma-\gamma^{a_{N-1}},$
$L_N:\gamma-\gamma^{a_{N-1}}\leq \gamma'<\gamma-1,$
$L_{N+1}:\gamma+1\leq\gamma'<\gamma+\gamma^{a_{N-1}},$
$L_{N+2}:\gamma+\gamma^{a_{N-1}}\leq\gamma'<\gamma+\gamma^{a_{N-2}},$
$\cdots,$
$L_{2N-1}:\gamma+\gamma^{a_2}\leq\gamma'<\gamma+\gamma^{a_1},$
$L_{2N}:\gamma+\gamma^{a_1}\leq \gamma'<2\gamma,$
$L_{2N+1}:2\gamma\leq\gamma'.$

 Note that, some of the $L_l$'s might be empty for those $\gamma$'s which are close to $T$, in which case the estimation will be trivial. Hence, we can assume each $L_l$ is not empty.

We take
$$a_1=1-\frac{1}{N},\quad a_2=1-\frac{2}{N},\quad \cdots,\quad a_{N-1}=\frac{1}{N}, \quad\mbox{and}\quad a_N=0.$$
Using the Cauchy-Schwarz inequality and (\ref{pf-lem-M-U-est}),
\begin{equation*}
  \sigma_1=\sum_{T\leq \gamma\leq X}\frac{|u_k(\rho)|}{|\rho\zeta'(\rho)|}\sum_{\gamma'\in L_1}\frac{|u_k(\rho')|}{|\rho'\zeta'(\rho')||\gamma-\gamma'|}
   \ll\sum_{T\leq \gamma\leq X}\frac{1}{|\gamma|^{\frac{1}{2}+\frac{1}{2k}+a_1-\epsilon}|\zeta'(\rho)|}\sum_{\gamma'\in L_1}\frac{1}{|\gamma'|^{\frac{1}{2}+\frac{1}{2k}-\epsilon}|\zeta'(\rho')|}.
\end{equation*}
Then, by partial summation and (\ref{J-1/2-est-pf-lem1}),
\begin{equation}\nonumber
 \sum_{\gamma'\in L_1}\frac{1}{|\gamma'|^{\frac{1}{2}+\frac{1}{2k}-\epsilon}|\zeta'(\rho')|}\ll \sum_{0<\gamma'<\gamma}\frac{1}{|\gamma'|^{\frac{1}{2}+\frac{1}{2k}-\epsilon}|\zeta'(\rho')|}\ll \gamma^{\frac{\theta}{2}-\frac{1}{2k}+2\epsilon}.
\end{equation}
Since $0<\epsilon<\frac{1+\frac{1}{k}-\theta}{10}$ and $N=[\frac{1}{\epsilon}]+2$,  $\frac{3}{2}+\frac{1}{k}-\frac{\theta}{2}-\frac{1}{N}-3\epsilon>\frac{\theta+1}{2}$, Lemma \ref{lem-G} applies. Thus, by (\ref{lem-M-assump}) and Lemma \ref{lem-G}, we obtain
\begin{equation}\label{sig-1}
 \sigma_1\ll \sum_{T\leq \gamma\leq X}\frac{\gamma^{\frac{\theta}{2}-\frac{1}{2k}+2\epsilon}}{|\gamma|^{\frac{1}{2}+\frac{1}{2k}+a_1-\epsilon}|\zeta'(\rho)|}
  \ll\sum_{T\leq \gamma\leq X}\frac{1}{|\gamma|^{\frac{3}{2}+\frac{1}{k}-\frac{\theta}{2}-\frac{1}{N}-3\epsilon}|\zeta'(\rho)|}
  \ll \frac{1}{T^{1+\frac{1}{k}-\theta-\frac{1}{N}-4\epsilon}}\ll \frac{1}{T^{1+\frac{1}{k}-\theta-5\epsilon}}.
\end{equation}

For $2\leq l\leq 2N$,  noting that $1+\frac{1}{k}-2\epsilon>\theta$, by the Cauchy-Schwarz inequality, (\ref{pf-lem-M-U-est}) and Lemma \ref{lem-G}, we have
\begin{eqnarray}
  \sigma_l&\ll& \left(\sum_{T\leq \gamma\leq X}\frac{|u_k(\rho)|^2}{|\rho\zeta'(\rho)|^2}\right)^{\frac{1}{2}}\left(\sum_{T\leq \gamma\leq X}\left(\sum_{\gamma'\in L_l}\frac{|u_k(\rho')|}{|\rho'\zeta'(\rho')||\gamma-\gamma'|}\right)^2\right)^{\frac{1}{2}}\nonumber\\
  &\ll&\left(\sum_{T\leq \gamma\leq X} \frac{1}{{\gamma}^{1+\frac{1}{k}-2\epsilon}|\zeta'(\rho)|^2}\right)^{\frac{1}{2}}\left(\sum_{T\leq \gamma\leq X}\left(\sum_{\gamma'\in L_l}\frac{|u_k(\rho')|^2}{|\rho'\zeta'(\rho')|^2|\gamma-\gamma'|^2}\right)N(L_l)\right)^{\frac{1}{2}}\nonumber\\
  &\ll&\frac{1}{T^{\frac{1}{2}+\frac{1}{2k}-\frac{\theta}{2}-\epsilon}}\left(\sum_{T\leq \gamma\leq X}\left(\sum_{\gamma'\in L_l}\frac{|u_k(\rho')|^2}{|\rho'\zeta'(\rho')|^2|\gamma-\gamma'|^2}\right)N(L_l)\right)^{\frac{1}{2}},\nonumber
\end{eqnarray}
where $N(L_l)$ is the number of $\gamma'$'s in $L_l$. By (\ref{zero-count}), $N(L_l)\ll \gamma^{a_{l-1}+\epsilon}$. Then, for $2\leq l\leq N$, since $\gamma'\asymp\gamma$ for $\gamma'\in L_l$, we have
\begin{equation}\nonumber
\frac{N(L_l)}{|\gamma-\gamma'|^2}\ll \frac{1}{\gamma^{2a_{l}-a_{l-1}-\epsilon}}=\frac{1}{\gamma^{1-\frac{l+1}{N}-\epsilon}}\ll\frac{1}{(\gamma')^{1-\frac{l+1}{N}-\epsilon}}.
\end{equation}
Then, we have
\begin{eqnarray}
  \sigma_l &\ll&\frac{1}{T^{\frac{1}{2}+\frac{1}{2k}-\frac{\theta}{2}-\epsilon}}\left(\sum_{T\leq \gamma\leq X}\sum_{\gamma'\in L_l}\frac{1}{|\gamma'|^{2+\frac{1}{k}-\frac{l+1}{N}-3\epsilon}|\zeta'(\rho')|^2}\right)^{\frac{1}{2}}\nonumber\\
  &\ll&\frac{1}{T^{\frac{1}{2}+\frac{1}{2k}-\frac{\theta}{2}-\epsilon}}\left(\sum_{T\leq \gamma'\leq X}m_l(\gamma')\frac{1}{|\gamma'|^{2+\frac{1}{k}-\frac{l+1}{N}-3\epsilon}|\zeta'(\rho')|^2}\right)^{\frac{1}{2}},\nonumber
\end{eqnarray}
by swapping summation and where $m_l(\gamma')=\#\{\gamma: \gamma-\gamma^{a_{l-1}}\leq \gamma'<\gamma-\gamma^{a_l}\}$. By (\ref{zero-count}),
$m_l(\gamma')=\#\{\gamma: \gamma'+\gamma^{a_{l}}< \gamma\leq\gamma'+\gamma^{a_{l-1}}\}
 \ll \#\{\gamma: \gamma'+(\gamma')^{a_l}<\gamma\leq \gamma'+(C\gamma')^{a_{l-1}} \text{~for some~} C>0 \}\ll(\gamma')^{1-\frac{l-1}{N}+\epsilon}$.

Since $N=[\frac{1}{\epsilon}]+2$ and $0<\epsilon<\frac{1+\frac{1}{k}-\theta}{10}$, $1+\frac{1}{k}-\frac{2}{N}-4\epsilon >\theta$, Lemma \ref{lem-G} applies. Thus, for $2\leq l\leq N$, we get,
\begin{equation}\label{sigma-l}
 \sigma_l\ll \frac{1}{T^{\frac{1}{2}+\frac{1}{2k}-\frac{\theta}{2}-\epsilon}}\left(\sum_{T\leq \gamma'\leq X}\frac{1}{|\gamma'|^{1+\frac{1}{k}-\frac{2}{N}-4\epsilon}|\zeta'(\rho')|^2}\right)^{\frac{1}{2}}\ll \frac{1}{T^{1+\frac{1}{k}-\theta-\frac{1}{N}-3\epsilon}}\ll\frac{1}{T^{1+\frac{1}{k}-\theta-4\epsilon}},
\end{equation}

Similarly, we have
\begin{equation}\label{sig-1-2}
  \sigma_{l}\asymp \sigma_{2N+1-l}, \text{~for~} N+1\leq l\leq 2N.
\end{equation}

Finally, we calculate $\sigma_{2N+1}$,
\begin{equation}\nonumber
  \sigma_{2N+1}\ll\sum_{T\leq \gamma\leq X}\frac{|u_k(\rho)|}{|\rho\zeta'(\rho)|}\left(\sum_{m=1}^{\infty}\sum_{2^m\gamma\leq\gamma'\leq 2^{m+1}\gamma}\frac{|u_k(\rho')|}{|\rho'\zeta'(\rho')||\gamma-\gamma'|}\right).
  \end{equation}
By (\ref{lem-M-assump}), (\ref{pf-lem-M-U-est}), and Lemma \ref{lem-G}, the inner sum is
\begin{eqnarray*}
  &\ll& \sum_{m=1}^{\infty}\frac{1}{(2^m-1)\gamma}\left(\sum_{2^m\gamma\leq\gamma'\leq 2^{m+1}\gamma}\frac{1}{|\rho'\zeta'(\rho')|^2}\right)^{\frac{1}{2}}\left(\sum_{2^m\gamma\leq\gamma'\leq 2^{m+1}\gamma}\left|u_k(\rho')\right|^2\right)^{\frac{1}{2}}\nonumber\\
  &\ll& \sum_{m=1}^{\infty}\frac{1}{(2^m-1)\gamma}\left(\frac{1}{2^m \gamma^{2-\theta}}\right)^{\frac{1}{2}}\left(2^{m+1}\gamma\right)^{\frac{1}{2}-\frac{1}{2k}+\epsilon}\left(2^{m+1}\gamma \log(2^{m+1}\gamma)\right)^{\frac{1}{2}}\nonumber\\
  &\ll& \sum_{m=1}^{\infty}\frac{1}{2^{(\frac{1}{2}+\frac{1}{2k}-\epsilon)m}}\frac{1}{\gamma^{1+\frac{1}{2k}-\frac{\theta}{2}-\epsilon}}
  \ll\frac{1}{\gamma^{1+\frac{1}{2k}-\frac{\theta}{2}-\epsilon}}.
\end{eqnarray*}
Noting that $\frac{3}{2}+\frac{1}{k}-\frac{\theta}{2}-2\epsilon>\frac{\theta+1}{2}$, by Lemma \ref{lem-G} and (\ref{pf-lem-M-U-est}),
\begin{equation}\label{sig-2N+1}
  \sigma_{2N+1}\ll \sum_{T\leq \gamma\leq X}\frac{|u_k(\rho)|}{|\rho\zeta'(\rho)|}\frac{1}{\gamma^{1+\frac{1}{2k}-\frac{\theta}{2}-\epsilon}} \ll \sum_{T\leq \gamma\leq X}\frac{1}{\gamma^{\frac{3}{2}+\frac{1}{k}-\frac{\theta}{2}-2\epsilon}|\zeta'(\rho)|}\ll \frac{1}{T^{1+\frac{1}{k}-\theta-3\epsilon}}.
\end{equation}

Combining all the estimates (\ref{sum2-sig}-\ref{sig-2N+1}),  we have

\begin{equation}\label{sum2-all}
  \Sigma_2\ll_{k, \epsilon} \frac{1}{T^{1+\frac{1}{k}-\theta-5\epsilon}}.
\end{equation}

Therefore, by (\ref{sum1}) and (\ref{sum2-all}), we deduce that, under the assumption (\ref{lem-M-assump}), for sufficiently small $\epsilon$,
\begin{equation}\nonumber
  \int_{\log Z}^{\log Z+1}\left|\sum_{T\leq \gamma\leq X}\frac{u_k(\rho)}{\rho\zeta'(\rho)}e^{\frac{i\gamma y}{k}}\right|^2 dy\ll_{k, \epsilon}\frac{1}{T^{1+\frac{1}{k}-\theta-\epsilon}}.
\end{equation}
The conclusion of Lemma \ref{lem-square} follows. \qed

\section{Proofs of Theorems}
\subsection{Proof of Theorem \ref{thm-average}: (\ref{GH}) implies (\ref{equiv-M}) }
By Lemma \ref{lemma-M}, for $X\leq x\ll X$, taking $T=X^2$, we have
\begin{equation}\nonumber
  M_k(x)=\sum_{|\gamma|<X^2}\frac{\zeta(\frac{\rho}{k})}{\rho\zeta'(\rho)}x^{\frac{\rho}{k}}+O(X^{\epsilon}).
\end{equation}
Then,
\begin{equation}\nonumber
  M^2_k(x)\ll \left|\sum_{|\gamma|<X^2}\frac{\zeta(\frac{\rho}{k})}{\rho\zeta'(\rho)}x^{\frac{\rho}{k}}\right|^2+O(X^{2\epsilon}).
\end{equation}

Since the imaginary part of the first zero of $\zeta(s)$ is $>14$, we let $x=e^y$, take $\theta=1+\epsilon$, $T=14$, $Z=X$, and replace $X$ by $X^2$ in Lemma \ref{lem-square}. We deduce that
\begin{equation}\label{form-thm1}
 \int_{X}^{eX} \left(\frac{M_k(x)}{x^{\frac{1}{2k}}}\right)^2 \frac{dx}{x}\ll \int_X^{eX}\left|\sum_{14<\gamma<X^2}\frac{\zeta(\frac{\rho}{k})}{\rho\zeta'(\rho)}x^{\frac{\rho}{k}}\right|^2 \frac{1}{x^{\frac{1}{k}}} \frac{dx}{x}+O(X^{-\frac{1}{k}+2\epsilon})\ll 1.
\end{equation}
So we get
\begin{equation}\nonumber
  \int_{2}^{X} \left(\frac{M_k(x)}{x^{\frac{1}{2k}}}\right)^2 \frac{dx}{x}\ll \sum_{l=1}^{[\log(\frac{X}{2})]+1}\int_{\frac{X}{e^{l-1}}}^{\frac{X}{e^l}} \left(\frac{M_k(x)}{x^{\frac{1}{2k}}}\right)^2 \frac{dx}{x}\ll_k \log X.
\end{equation}

Here, it is easy to get a weaker result than Theorem \ref{thm-asymp}. By (\ref{form-thm1}), we have
\begin{equation}\nonumber
  \int_{X}^{eX} \left(\frac{M_k(x)}{x^{\frac{1}{2k}}}\right)^2 dx\ll X.
\end{equation}
Substituting $\frac{X}{e},~\frac{X}{e^2},~\cdots$, for $X$ in the above formula, we obtain
\begin{equation}\nonumber
  \int_{2}^{X} \left(\frac{M_k(x)}{x^{\frac{1}{2k}}}\right)^2 dx\ll X.
\end{equation}

\subsection{Proof of Theorem \ref{thm-average}: (\ref{equiv-M}) implies (\ref{GH})}
We need the following lemmas.
\begin{lem}\label{equiv-lem-O}
Assume the Riemann Hypothesis. The formula (\ref{equiv-M}) implies the zeros of $\zeta(s)$ on the critical line are simple, that
\begin{equation}\nonumber
  \frac{\zeta(\frac{\rho}{k})}{\zeta'(\rho)}=O_k(|\rho|), \quad\forall k\geq 2,
\end{equation}
and
\begin{equation}\nonumber
 \frac{1}{\zeta'(\rho)}=O_{\epsilon}\left(|\rho|^{\frac{1}{2}+\epsilon}\right), \quad \forall \epsilon>0.
\end{equation}
\end{lem}

 Using the above lemma and a similar argument to Theorem 14.29 (B) of \cite{Titch}, we prove the following Lemma.
\begin{lem}\label{equiv-lem-conv}
Assume the Riemann Hypothesis. If (\ref{equiv-M}) is true, then for any $k\geq 2$ the series
  \begin{equation*}
    \sum_{\rho}\frac{|\zeta(\frac{\rho}{k})|^2}{|\rho\zeta'(\rho)|^2}
  \end{equation*}
  is convergent.
\end{lem}

Before giving the proofs of these two lemmas, we show how Lemma \ref{equiv-lem-conv} implies the desired result.

For any $\epsilon>0$, under the Riemann Hypothesis, we know that (\cite{Titch}, (14.2.6), p.337)
\begin{equation}\nonumber
  \frac{1}{\zeta(\sigma+it)}=O(|t|^{\frac{\epsilon}{4}}), \mbox{~for every fixed~} \sigma>\frac{1}{2}.
\end{equation}
So by the functional equation, we have
\begin{equation}\label{pf-thm-equiv-conv-zeta}
  |\zeta(\frac{\rho}{k})|\gg_k |\rho|^{\frac{1}{2}-\frac{1}{2k}-\frac{\epsilon}{4}}.
\end{equation}
Taking $k=\left[\frac{2}{\epsilon}\right]+2$, then
\begin{equation}\nonumber
  |\zeta(\frac{\rho}{k})|\gg_{\epsilon} |\rho|^{\frac{1}{2}-\frac{\epsilon}{2}}.
\end{equation}
Then, by Lemma \ref{equiv-lem-conv}, we have
\begin{equation}\nonumber
  \sum_{0<\gamma\leq T}\frac{1}{\gamma^{1+\epsilon}|\zeta'(\rho)|^2}\ll \sum_{0<\gamma\leq T}\frac{|\zeta(\frac{\rho}{k})|^2}{|\rho\zeta'(\rho)|^2}\ll_{\epsilon} 1.
\end{equation}

Hence, by partial summation, we derive that, for any $\epsilon>0$,
\begin{equation}\nonumber
  \sum_{0<\gamma\leq T}\frac{1}{|\zeta'(\rho)|^2}=\sum_{0<\gamma\leq T}\frac{1}{\gamma^{1+\epsilon}|\zeta'(\rho)|^2}\gamma^{1+\epsilon}\ll_{\epsilon} T^{1+\epsilon}+\int_1^T t^{\epsilon}dt\ll_{\epsilon} T^{1+\epsilon}.
\end{equation}
This proves (\ref{GH}). \qed

\vspace{1em}

\noindent \emph{Proof of Lemma \ref{equiv-lem-O}.} By considering Mellin transforms, for $\Re s>1$,
\begin{equation}\label{pf-thm-equiv-lem1-mellin}
  \frac{\zeta(s)}{\zeta(ks)}-\frac{\zeta(s)}{\zeta(k)}=\sum_{n=1}^{\infty}\frac{\mu_k(n)-\frac{1}{\zeta(k)}}{n^s}  =s\int_1^{\infty}\frac{M_k(x)+\frac{\{x\}}{\zeta(k)}}{x^{s+1}}dx
  =s\int_1^{\infty}\frac{M_k(x)}{x^{s+1}}dx+\frac{s}{\zeta(k)}\int_1^{\infty}\frac{\{x\}}{x^{s+1}}dx.
\end{equation}
We show that the right hand side of (\ref{pf-thm-equiv-lem1-mellin}) can be analytically continued to the half-plane $\Re s>\frac{1}{2k}$. Let $s=\sigma+it$ with $\sigma>\frac{1}{2k}$, we have
\begin{eqnarray}\label{pf-thm-equiv-lem1-cov}
  \left|s\int_1^{\infty}\frac{M_k(x)}{x^{s+1}}dx\right|&\leq&|s|\int_1^{\infty}\frac{|M_k(x)|}{x^{\sigma+1}}dx
  =|s|\int_1^{\infty}\frac{|M_k(x)|}{x^{\frac{1}{2}\sigma+\frac{1}{2}+\frac{1}{4k}}}\frac{1}{x^{\frac{1}{2}\sigma+\frac{1}{2}-\frac{1}{4k}}}dx\nonumber\\
  &\leq& |s|\left(\int_1^{\infty}\frac{M^2_k(x)}{x^{\sigma+1+\frac{1}{2k}}}dx\right)^{\frac{1}{2}}\left(\int_1^{\infty}\frac{1}{x^{\sigma+1-\frac{1}{2k}}}dx\right)^{\frac{1}{2}}\nonumber\\
  &\leq&\frac{|s|}{\sqrt{\sigma-\frac{1}{2k}}}\left(\int_1^{\infty}\frac{M^2_k(x)}{x^{\sigma+1+\frac{1}{2k}}}dx\right)^{\frac{1}{2}}.
\end{eqnarray}
Let
\begin{equation}\nonumber
  f(X)=\int_1^X \frac{M^2_k(x)}{x^{1+\frac{1}{k}}}dx.
\end{equation}
By assumption (\ref{equiv-M}), $f(X)\ll \log X$. Then using integration by parts,
\begin{eqnarray}\label{h-Mk-est}
  \int_1^{\infty}\frac{M^2_k(x)}{x^{\sigma+1+\frac{1}{2k}}}dx&=&\int_1^{\infty}\frac{f'(x)}{x^{\sigma-\frac{1}{2k}}}dx
  =\left(\sigma-\frac{1}{2k}\right)\int_1^{\infty}\frac{f(x)}{x^{\sigma+1-\frac{1}{2k}}}dx\nonumber\\
  &=&O\left(\left(\sigma-\frac{1}{2k}\right)\int_1^{\infty}\frac{\log x}{x^{\sigma+1-\frac{1}{2k}}}dx\right)
  =O\left(\int_1^{\infty}\frac{1}{x^{\sigma+1-\frac{1}{2k}}}dx\right)\nonumber\\
  &=&O\left(\frac{1}{\sigma-\frac{1}{2k}}\right).
\end{eqnarray}
Then, by (\ref{pf-thm-equiv-lem1-cov}) and (\ref{h-Mk-est}), we get
\begin{equation}\nonumber
  \left|s\int_1^{\infty}\frac{M_k(x)}{x^{s+1}}dx\right|=O\left(\frac{|s|}{\sigma-\frac{1}{2k}}\right),
\end{equation}
which shows that the integral converges uniformly and absolutely for $\sigma\geq \frac{1}{2k}+\delta$, for any fixed $\delta>0$. The second integral on the right hand side of (\ref{pf-thm-equiv-lem1-mellin}) is also uniformly and absolutely convergent for $\sigma\geq \frac{1}{2k}$, and
\begin{equation}\nonumber
 \left| \frac{s}{\zeta(k)}\int_1^{\infty}\frac{\{x\}}{x^{s+1}}dx\right|\ll 1.
\end{equation}
Thus, the formula (\ref{pf-thm-equiv-lem1-mellin}) can be analytically continued to $\Re s>\frac{1}{2k}$, and we get
\begin{equation}\label{ratio-zeta-est2}
 \frac{\zeta(s)}{\zeta(ks)}-\frac{\zeta(s)}{\zeta(k)}=O\left(\frac{|s|}{\sigma-\frac{1}{2k}}\right)+O(1).
\end{equation}

Let $\rho$ be a nontrivial zero of $\zeta(s)$, and $s=\frac{\rho}{k}+\frac{h}{k}$, $h>0$ be small real positive number. Then by (\ref{ratio-zeta-est2}),
\begin{equation}\label{ratio-est-h}
  \frac{\zeta(\frac{\rho}{k}+\frac{h}{k})}{\zeta(\rho+h)}=O\left(\frac{|\frac{\rho}{k}+\frac{h}{k}|}{\frac{h}{k}}\right)
  +O\left(\frac{|\zeta(\frac{\rho}{k}+\frac{h}{k})|}{|\zeta(k)|}+1\right)
  =O\left(\frac{|\rho+h|}{h}\right)+O\left(\frac{|\zeta(\frac{\rho}{k}+\frac{h}{k})|}{|\zeta(k)|}+1\right).
\end{equation}
This would be false for $h\rightarrow 0$ if $\rho$ is not a simple zero. Multiplying by $h$ on both sides of (\ref{ratio-est-h}), and letting $h\rightarrow 0$, we get
\begin{equation}\nonumber
 \frac{\zeta(\frac{\rho}{k})}{\zeta'(\rho)}=O(|\rho|),
\end{equation}
where the constant in big $O$ depends on $k$.

For any $\epsilon>0$, by (\ref{pf-thm-equiv-conv-zeta}),
$\zeta(\frac{\rho}{k})\gg |\rho|^{\frac{1}{2}-\frac{1}{2k}-\frac{\epsilon}{2}}.$
Then, we have
  $\frac{1}{\zeta'(\rho)}=O\left(|\rho|^{\frac{1}{2}+\frac{1}{2k}+\frac{\epsilon}{2}}\right)$.
Taking $k=\left[\frac{1}{\epsilon}\right]+2$, hence we get
$\frac{1}{\zeta'(\rho)}=O\left(|\rho|^{\frac{1}{2}+\epsilon}\right)$. \qed

\vspace{1em}
\noindent\emph{Proof of Lemma \ref{equiv-lem-conv}.} By the symmetry of $\zeta(s)$, the sum in the following formula is actually real. We have
\begin{eqnarray}\label{lem-cov-main-ineq}
  0&\leq&\int_1^X\left(\frac{M_k(x)}{x^{\frac{1}{2k}}}-\sum_{|\gamma|<T}\frac{\zeta(\frac{\rho}{k})x^{\frac{\rho}{k}-\frac{1}{2k}}}{\rho\zeta'(\rho)}\right)^2 \frac{dx}{x}\nonumber\\
  &=&\int_1^X \left(\frac{M_k(x)}{x^{\frac{1}{2k}}}\right)^2\frac{dx}{x}+\sum_{|\gamma|,|\gamma'|<T}\frac{\zeta(\frac{\rho}{k})\zeta(\frac{\rho'}{k})}{\rho\rho'\zeta'(\rho)\zeta'(\rho')}\int_1^X x^{\frac{\rho+\rho'-1}{k}}\frac{dx}{x}\nonumber\\
  &&-2\sum_{|\gamma|<T}\frac{\zeta(\frac{\rho}{k})}{\rho\zeta'(\rho)}\int_1^X M_k(x)x^{\frac{\rho-1}{k}}\frac{dx}{x}.
\end{eqnarray}

In the first sum of (\ref{lem-cov-main-ineq}), the terms with $\rho'=1-\rho$ contribute
\begin{equation}
  \sum_{|\gamma|<T}\frac{\zeta(\frac{\rho}{k})\zeta(\frac{1-\rho}{k})}{\rho(1-\rho)\zeta'(\rho)\zeta'(1-\rho)}\int_1^X \frac{dx}{x}=\log X \sum_{|\gamma|<T}\frac{|\zeta(\frac{\rho}{k})|^2}{|\rho\zeta'(\rho)|^2}.
\end{equation}
For the remaining terms, we write $\rho=\frac{1}{2}+i\gamma$, $\rho'=\frac{1}{2}+i\gamma'$, and $\gamma'\neq -\gamma$,
\begin{equation}
  \int_1^X x^{\frac{\rho+\rho'-1}{k}}\frac{dx}{x}=\frac{k(X^{\frac{\rho+\rho'-1}{k}}-1)}{\rho+\rho'-1}=O\left(\frac{k}{|\gamma+\gamma'|}\right).
\end{equation}
Thus, the sum of these terms is less than a constant $K_1=K_1(k, T)$.

In the last sum of (\ref{lem-cov-main-ineq}),
\begin{equation}\label{M-split-est}
  \int_1^X M_k(x) x^{\frac{\rho-1}{k}}\frac{dx}{x}=\int_1^X M_k(x)x^{\frac{\rho-1}{k}}\left(1-\frac{x}{X}\right)\frac{dx}{x}+\frac{1}{X}\int_1^X M_k(x)x^{\frac{\rho-1}{k}}dx.
\end{equation}
By the assumption (\ref{equiv-M}), the last term is
\begin{eqnarray}\label{M-split-est-1st}
  \frac{1}{X}\int_1^X M_k(x)x^{\frac{\rho-1}{k}}dx&=&O\left(\frac{1}{X}\int_1^X |M_k(x)|x^{-\frac{1}{2k}}dx\right)
  =O\left(\frac{1}{X}\int_1^X \frac{|M_k(x)|}{x^{\frac{1}{2k}}}\frac{1}{\sqrt{x}}\sqrt{x}dx\right)\nonumber\\
  &=&O\left(\frac{1}{X}\left(\int_1^X \left(\frac{M_k(x)}{x^{\frac{1}{2k}}}\right)^2\frac{dx}{x}\right)^{\frac{1}{2}}\left(\int_1^X xdx\right)^{\frac{1}{2}}\right)
  =O(\sqrt{\log X}).
\end{eqnarray}
For the first term in (\ref{M-split-est}),
\begin{eqnarray}\label{M-split-2nd-est}
  &&\int_1^X M_k(x)x^{\frac{\rho-1}{k}}\left(1-\frac{x}{X}\right)\frac{dx}{x}
  =\int_1^X \left(\sum_{n\leq x}\mu_k(n)-\frac{x}{\zeta(k)}\right)x^{\frac{\rho-1}{k}}\left(1-\frac{x}{X}\right)\frac{dx}{x}\nonumber\\
  &&\qquad=\int_1^X \left(\sum_{n\leq x}\mu_k(n)\right)x^{\frac{\rho-1}{k}}\left(1-\frac{x}{X}\right)\frac{dx}{x}-\frac{1}{\zeta(k)}\int_1^X\left(1-\frac{x}{X}\right)x^{\frac{\rho-1}{k}}dx.
\end{eqnarray}
For the last integral,
\begin{eqnarray}
  \frac{1}{\zeta(k)}\int_1^X\left(1-\frac{x}{X}\right)x^{\frac{\rho-1}{k}}dx
  =\frac{k^2}{\zeta(k)}\frac{X^{\frac{\rho-1}{k}+1}-1}{(\rho-1+k)(\rho-1+2k)}-\frac{k}{\zeta(k)(\rho-1+2k)}\left(1-\frac{1}{X}\right).
\end{eqnarray}
For the first integral in (\ref{M-split-2nd-est}), we show that
\begin{equation}\label{M-split-2nd2int-est}
  \int_1^X \left(\sum_{n\leq x}\mu_k(n)\right)x^{\frac{\rho-1}{k}}\left(1-\frac{x}{X}\right)\frac{dx}{x}=\frac{1}{2\pi i}\int_{2k-i\infty}^{2k+i\infty}k^2\frac{\zeta(\frac{w}{k})}{\zeta(w)}\frac{X^{\frac{w+\rho-1}{k}}-1}{w(w+\rho-1+k)(w+\rho-1)}dw.
\end{equation}
To prove this, we use the following formula, for $a>\frac{1}{2}$, $y>0$,
\begin{equation}\label{my-perron-for}
  \frac{1}{2\pi i}\int_{a-i\infty}^{a+i\infty} \frac{y^s}{s(s+\rho-1+k)(s+\rho-1)} ds
=\left\{
    \begin{array}{ll}                                                                                      \frac{1}{k}\left(\frac{y^{1-\rho}-1}{1-\rho}+\frac{y^{1-\rho-k}-1}{\rho-1+k}\right), & \hbox{~for~} y>1;\\
                                                                                         0 , & \hbox{~for~}0<y\leq 1.
                                                                                        \end{array}
                                                                                      \right.
\end{equation}
In fact, for $y>1$, we shift the path of the integration to the left as far as possible and calculate the residues at $s=0$, $s=1-\rho$, and $s=1-\rho-k$. Then we get the first result. For $0<y<1$, we shift the path of the integration to the right as far as possible. Since there are no poles to the right of the line $\Re s=a>\frac{1}{2}$ and $0<y<1$, the integral is zero. For $y=1$, by the continuity of $y$ on both sides of (\ref{my-perron-for}), we see that the integral equals zero.

By the abosolute convergence of $ \frac{\zeta(\frac{w}{k})}{\zeta(w)}=\sum_{n=1}^{\infty}\frac{\mu_k(n)}{n^{\frac{w}{k}}}$, we substitute it
on the right hand side of (\ref{M-split-2nd2int-est}) and integrate term by term.
By (\ref{my-perron-for}), taking $y=\left(\frac{X}{n}\right)^{\frac{1}{k}}$,
we obtain
\begin{eqnarray}
  &&\sum_{n=1}^{\infty}\frac{\mu_k(n)}{2\pi i}\int_{2k-i\infty}^{2k+i\infty}\frac{k^2 }{n^{\frac{w}{k}}}\frac{X^{\frac{w+\rho-1}{k}}}{w(w+\rho-1+k)(w+\rho-1)}dw\nonumber\\
  &&\qquad=\sum_{n\leq X}\mu_k(n)\left(k\frac{X^{\frac{\rho-1}{k}}-n^{\frac{\rho-1}{k}}}{\rho-1}-\frac{k}{\rho-1+k}\frac{X^{\frac{\rho-1}{k}+1}-n^{\frac{\rho-1}{k}+1}}{X}\right)\nonumber\\
  &&\qquad=\sum_{n\leq X}\mu_k(n)\int_n^X\left(\frac{x^{\frac{\rho-1}{k}}}{x}-\frac{x^{\frac{\rho-1}{k}}}{X}\right)dx\nonumber\\
  &&\qquad=\int_1^X \left(\sum_{n\leq x}\mu_k(n)\right)x^{\frac{\rho-1}{k}}\left(1-\frac{x}{X}\right)\frac{dx}{x}.\nonumber
\end{eqnarray}
And taking $y=\left(\frac{1}{n}\right)^{\frac{1}{k}}$, we get, for all $n\geq 1$,
$$\int_{2k-i\infty}^{2k+i\infty}\frac{k^2 }{n^{\frac{w}{k}}}\frac{1}{w(w+\rho-1+k)(w+\rho-1)}dw=0.$$
Thus,
$$\sum_{n=1}^{\infty}\frac{\mu_k(n)}{2\pi i}\int_{2k-i\infty}^{2k+i\infty}\frac{k^2 }{n^{\frac{w}{k}}}\frac{X^{\frac{w+\rho-1}{k}}-1}{w(w+\rho-1+k)(w+\rho-1)}dw
=\int_1^X \left(\sum_{n\leq x}\mu_k(n)\right)x^{\frac{\rho-1}{k}}\left(1-\frac{x}{X}\right)\frac{dx}{x},$$
which is the left hand side of (\ref{M-split-2nd2int-est}).

Let $U>T$ such that U is not the ordinate of a zero of $\zeta(s)$. Then the right hand side of (\ref{M-split-2nd2int-est}) is equal to
\begin{equation}\label{RHS-sum2int}
  \frac{1}{2\pi i}\left(\int_{2k-i\infty}^{2k-iU}+\int_{2k-iU}^{\frac{1}{4}-iU}+\int_{\frac{1}{4}-iU}^{\frac{1}{4}+iU}+\int_{\frac{1}{4}+iU}^{2k+iU}+\int_{2k+iU}^{2k+i\infty}\right)+\mbox{sum of residues in} -U<\Im w<U.
\end{equation}

Let $\rho''$ be a generic zero of $\zeta(s)$ with $|\gamma''|<U$. Since $U>T$, $w=1-\rho$ is a pole with residue
\begin{equation}
  \frac{\zeta(\frac{1-\rho}{k})}{(1-\rho)\zeta'(1-\rho)}\log X.
\end{equation}
There is also a pole at $w=k$ with residue
\begin{equation}
  \frac{k^2}{\zeta(k)}\frac{X^{\frac{\rho-1}{k}+1}-1}{(\rho-1+k)(\rho-1+2k)}.
\end{equation}
The residue at other $\rho''$ is
\begin{equation}\nonumber
  \mbox{Res} (\rho'')=k^2\frac{\zeta(\frac{\rho''}{k})}{\zeta'(\rho'')}\frac{X^{\frac{\rho''+\rho-1}{k}}-1}{\rho''(\rho''+\rho-1+k)(\rho''+\rho-1)}.
\end{equation}
By Lemma \ref{equiv-lem-O},
\begin{equation}\nonumber
  \frac{\zeta(\frac{\rho''}{k})}{\zeta'(\rho'')}=O_k(|\rho''|).
\end{equation}
Then,
\begin{equation}
  \mbox{Res} (\rho'')=O_k\left(\frac{1}{|(\rho''+\rho-1+k)(\rho''+\rho-1)|}\right)=O_k\left(\frac{1}{|\gamma+\gamma''|^2}\right).
\end{equation}
Thus, since $|\gamma|<T$,
\begin{equation}
  \sum_{\substack{-U<\gamma''<U\\ \gamma''\neq -\gamma}}\frac{1}{|\gamma''+\gamma|^2}\leq \sum_{\gamma''\neq-\gamma}\frac{1}{|\gamma''+\gamma|^2}<K_2(T).
\end{equation}

In the following, we estimate those five integrals in (\ref{RHS-sum2int}).
First, we have
\begin{eqnarray}
  &&\int_{2k+iU}^{2k+i\infty}k^2\frac{\zeta(\frac{w}{k})}{\zeta(w)}\frac{X^{\frac{w+\rho-1}{k}}-1}{w(w+\rho-1+k)(w+\rho-1)}dw\nonumber\\
  &&\quad=O_k\left(X^2\int_U^{\infty}\frac{dv}{v(v+\gamma)^2}\right)=O_k\left(\frac{X^2}{U(U+\gamma)}\right)=O_k\left(\frac{X^2}{U(U-T)}\right).
\end{eqnarray}
Similarly, we get the same estimate for the integral over $(2k-i\infty, 2k-iU)$. Next,
\begin{eqnarray}
  &&\int_{\frac{1}{4}-iU}^{\frac{1}{4}+iU}k^2\frac{\zeta(\frac{w}{k})}{\zeta(w)}\frac{X^{\frac{w+\rho-1}{k}}-1}{w(w+\rho-1+k)(w+\rho-1)}dw\nonumber\\
  &&\quad =\left(\int_{\frac{1}{4}-iT}^{\frac{1}{4}+iT}+\int_{\frac{1}{4}-iU}^{\frac{1}{4}-iT}+\int_{\frac{1}{4}+iT}^{\frac{1}{4}+iU}\right)
  k^2\frac{\zeta(\frac{w}{k})}{\zeta(w)}\frac{X^{\frac{w+\rho-1}{k}}-1}{w(w+\rho-1+k)(w+\rho-1)}dw\nonumber\\
  &&\quad =K'_3(k, T)+O\left(\int_T^U \frac{dv}{v^{\frac{1}{2}}(v+\gamma)^2}\right)
  =K'_3(k, T)+O(1)
   \leq K_3(k, T).
\end{eqnarray}
By Lemma \ref{lemma3-Ng}, we choose $n\leq U=U_n\leq n+1$ so that
\begin{equation}\nonumber
  \frac{1}{\zeta(\sigma+iU)}=O(|U|^{\epsilon}), \mbox{~for~}\frac{1}{4}\leq \sigma\leq 2.
\end{equation}
and we have
$\zeta(\frac{\sigma+iU}{k})=O(U^{\frac{1}{2}+\epsilon})$.
Then, we get
\begin{equation}\label{last-int-est}
  \int_{\frac{1}{4}+iU}^{2k+iU}k^2\frac{ \zeta(\frac{w}{k})}{\zeta(w)}\frac{X^{\frac{w+\rho-1}{k}}-1}{w(w+\rho-1+k)(w+\rho-1)}dw
 =O_k\left(\frac{X^2}{U^{\frac{1}{2}-\epsilon}(U+\gamma)^2}\right)=O_k\left(\frac{X^2}{U^{\frac{1}{2}-\epsilon}(U-T)^2}\right).
\end{equation}
Similarly for the integral over $(2k-iU, \frac{1}{4}-iU)$.

Combining (\ref{M-split-2nd-est}-\ref{last-int-est}), and making $U\rightarrow\infty$, we get
\begin{equation}\label{M-split-2nd-res}
  \int_1^X M_k(x)x^{\frac{\rho-1}{k}}\left(1-\frac{x}{X}\right)\frac{dx}{x}=\frac{\zeta(\frac{1-\rho}{k})}{(1-\rho)\zeta'(1-\rho)}\log X+R,
\end{equation}
where $|R|<K_4(k, T)$ if $|\gamma|<T$.

By (\ref{lem-cov-main-ineq}-\ref{M-split-est-1st}) and (\ref{M-split-2nd-res}), and by assumption (\ref{equiv-M}), we deduce that
\begin{equation}\nonumber
  0\leq A_k \log X+\log X\sum_{|\gamma|<T}\frac{|\zeta(\frac{\rho}{k})|^2}{|\rho\zeta'(\rho)|^2}-2\log X\sum_{|\gamma|<T}\frac{|\zeta(\frac{\rho}{k})|^2}{|\rho\zeta'(\rho)|^2}+A'_k\sqrt{\log X}+K(k, T),
\end{equation}
where $A_k$ and $A'_k$ are constants depending only on $k$.
Thus,
\begin{equation}\nonumber
  \sum_{|\gamma|<T}\frac{|\zeta(\frac{\rho}{k})|^2}{|\rho\zeta'(\rho)|^2}\leq A_k+\frac{A'_k}{\sqrt{\log X}}+\frac{K(k, T)}{\log X}.
\end{equation}
Making $X\rightarrow\infty$,
$\sum_{|\gamma|<T}\frac{|\zeta(\frac{\rho}{k})|^2}{|\rho\zeta'(\rho)|^2}\leq A_k$.
Since the right hand side of the above result is independent of $T$, we get the convergence of
$\sum_{\rho}\frac{|\zeta(\frac{\rho}{k})|^2}{|\rho\zeta'(\rho)|^2}$. \qed

\subsection{Proof of Theorem \ref{thm-asymp}}
We define
$\phi(y)=e^{-\frac{y}{2k}}M_k(e^{y})=\phi^{(T)}(y)+\epsilon^{(T)}(y)$,
where
\begin{equation}\nonumber
  \phi^{(T)}(y)=\sum_{|\gamma|\leq T}\frac{\zeta(\frac{\rho}{k})}{\rho\zeta'(\rho)}e^{\frac{i\gamma y}{k}},
\end{equation}
and, letting $Y=\log X$,
\begin{equation}\label{pf-thm3-epsi}
  \epsilon^{(T)}(y)=\sum_{T<|\gamma|\leq e^{2Y}}\frac{\zeta(\frac{\rho}{k})}{\rho\zeta'(\rho)}e^{\frac{i\gamma y}{k}}+e^{-\frac{y}{2k}}E(e^{y}, e^{2Y}),
\end{equation}
where $E(x, T)$ is as defined in Lemma \ref{lemma-M}.
Note that
\begin{equation}\nonumber
\int_{\log 2}^{Y}|e^{-\frac{y}{2k}}E(e^{y}, e^{2Y})|^2dy\ll\int_{\log 2}^{Y}\frac{y^2 e^{(2-\frac{1}{k})y}}{e^{4Y}}+\frac{e^{(2-\frac{1}{k})y}}{e^{2(1-\epsilon)Y}}+\frac{1}{e^{4(\frac{1}{2k}-\epsilon)Y}}+\frac{1}{e^{\frac{y}{k}}} dy\ll 1.
\end{equation}
For $k\geq 2$, taking $\theta=1+\epsilon$ in Lemma \ref{lem-square}, and by (\ref{pf-thm3-epsi}), we have
\begin{eqnarray*}
  &&\limsup_{Y\rightarrow\infty}\frac{1}{Y}\int_{0}^Y \left|\phi(y)-\sum_{|\gamma|\leq T}\frac{\zeta(\frac{\rho}{k})}{\rho\zeta'(\rho)}e^{\frac{i\gamma y}{k}}\right|^2 dy=	\limsup_{Y\rightarrow\infty}\frac{1}{Y}\left(\int_{0 }^{\log 2}+\int_{\log 2}^{Y}\right)|\epsilon^{(T)}(y)|^2dy\\
  &&\qquad \ll\limsup_{Y\rightarrow\infty}\frac{1}{Y}\int_{\log 2}^{Y}\left|\sum_{T\leq \gamma\leq e^{2Y}}\frac{\zeta(\frac{\rho}{k})}{\rho\zeta'(\rho)}e^{\frac{i\gamma y}{k}}\right|^2 dy+\lim_{Y\rightarrow\infty}\frac{1}{Y}\int_{\log 2}^{Y}|e^{-\frac{y}{2k}}E(e^{y}, e^{2Y})|^2dy\\
  &&\qquad\ll \limsup_{Y\rightarrow\infty}\frac{1}{Y}\sum_{j=0}^{[Y]}\int_{\log 2+j}^{\log 2+j+1}\left|\sum_{T\leq\gamma\leq e^{2Y}}\frac{\zeta(\frac{\rho}{k})}{\rho\zeta'(\rho)}e^{\frac{i\gamma y}{k}}\right|^2dy\ll \frac{1}{T^{\frac{1}{k}-\epsilon}}.\nonumber
\end{eqnarray*}

Then, by (\ref{B2-almost}), we see that $\phi(y)=e^{-\frac{y}{2k}}M_k(e^{y})$ is a $B^2$-almost periodic function. Thus, by the work of Besicovitch (see \cite{Besi}, Chapter II of \cite{Besi2}, or Theorem 1.14 of \cite{A-N-S}), we get the conclusion of our theorem. \qed

\section{Applications of LI}
In this section, we assume the Linear Independence conjecture, and give the proof of Corollary \ref{app-LI}.

Let $X$ be a random variable on the infinite torus $\mathbb{T}^{\infty}$,
$$X(\boldsymbol{\theta})=\sum_{l=1}^{\infty}r_l \sin 2\pi\theta_l,$$
where $\boldsymbol{\theta}=(\theta_1, \theta_2, \cdots)\in \mathbb{T}^{\infty}$ and $r_l\in \mathbb{R}$ for $l\geq 1$. If we assume $\sum_{l=1}^{\infty}r^2_l<\infty$, then $X$ converges almost everywhere by Komolgorov's theorem.
Let $\mathbf{P}$ be the canonical probability measure on $\mathbb{T}^{\infty}$. Define
$$\nu_X(x)=\mathbf{P}(X^{-1}(-\infty, x)).$$

Montgomery \cite{Mont} proved the following result. It also follows from Hoeffding's inequality \cite{Hoeff}.
\begin{lem}\label{lem-montg}
  Let $X(\boldsymbol{\theta})=\sum_{l=1}^{\infty}r_l \sin 2\pi\theta_l$ where $\sum_{l=1}^{\infty}r^2_l<\infty$. For any integer $K\geq 1$,
  \begin{equation*}
    \mathbf{P}\left(X(\boldsymbol{\theta})\geq 2\sum_{l=1}^{K}r_l\right)\leq \exp\left(-\frac{3}{4}\left(\sum_{l=1}^{K}r_l\right)^2\left(\sum_{l>K} r_l^2\right)^{-1}\right).
  \end{equation*}
\end{lem}

The linear independence assumption implies that the limiting distribution $\nu$ obtained in Theorem \ref{thm-distri} equals $\nu_X$, where $X$ is the random variable
\begin{equation}\nonumber
  X(\boldsymbol{\theta})=\sum_{\gamma>0}\frac{2|\zeta(\frac{\rho}{k})|}{|\rho\zeta'(\rho)|}\sin(2\pi\theta_{\gamma}).
\end{equation}
Let $r_{\gamma}=\frac{2|\zeta(\frac{\rho}{k})|}{|\rho\zeta'(\rho)|}$. Define
\begin{equation}\nonumber
  A(T):=\sum_{0<\gamma<T}r_{\gamma}=\sum_{0<\gamma<T}\frac{2|\zeta(\frac{\rho}{k})|}{|\rho\zeta'(\rho)|}, \quad \mbox{and}\quad B(T):=\sum_{\gamma\geq T}r^2_{\gamma}=\sum_{\gamma\geq T}\frac{4|\zeta(\frac{\rho}{k})|^2}{|\rho\zeta'(\rho)|^2}.
\end{equation}
By $(1. 134)$ in \cite{Ivic} (pp. 45), $(14.2.6)$ and $(14.5.1)$ in \cite{Titch}, and the functional equation,
\begin{equation}\nonumber
  \left|\zeta(\frac{\rho}{k})\right|\gg |\gamma|^{\frac{1}{2}-\frac{1}{2k}-\epsilon}, \quad\mbox{and}\quad |\zeta'(\rho)|\ll |\gamma|^{\epsilon}.
\end{equation}
So, we deduce that, for $\gamma<T$,
\begin{equation}\label{r-gamma-est}
  r_{\gamma}=\frac{2|\zeta(\frac{\rho}{k})|}{|\rho\zeta'(\rho)|}\gg_k \frac{1}{|\gamma|^{\frac{1}{2}+\frac{1}{2k}+\epsilon}}\gg \frac{1}{T^{\frac{1}{2}+\frac{1}{2k}+\epsilon}}.
\end{equation}
Then, by (\ref{zero-count}),
\begin{equation}\nonumber
  A(T)\gg \sum_{0<\gamma<T}\frac{1}{T^{\frac{1}{2}+\frac{1}{2k}+\epsilon}}\gg T^{\frac{1}{2}-\frac{1}{2k}-\epsilon}.
\end{equation}
Thus, by partial summation, (\ref{GH}),  and the Riemann Hypothesis, we get
\begin{equation}\label{conj-formula}
T^{\frac{1}{2}-\frac{1}{2k}-\epsilon}\ll A(T)\ll T^{\frac{1}{2}-\frac{1}{2k}+\epsilon}, \quad\mbox{and}\quad \frac{1}{T^{\frac{1}{k}+\epsilon}}\ll B(T)\ll \frac{1}{T^{\frac{1}{k}-\epsilon}}.
\end{equation}

Let $V$ be a large parameter. We want to use the above estimates to find the upper and lower bounds for the tail of the distribution,
\begin{equation}\nonumber
  \nu([V, \infty)):=\int_{V}^{\infty}d\nu(x)=\mathbf{P}(X(\boldsymbol{\theta})\geq V).
\end{equation}

\subsection{The upper bound}
Choose $T$ such that $A(T^-)<V\leq A(T)$. We have the inequalities,
\begin{equation}\nonumber
T^{\frac{1}{2}-\frac{1}{2k}-\epsilon}\ll A(T^-)<V\leq A(T)\ll T^{\frac{1}{2}-\frac{1}{2k}+\epsilon}.
\end{equation}
From this, we see that
\begin{equation}\label{T-V}
V^{\frac{2k}{k-1}-\epsilon}\ll T\ll V^{\frac{2k}{k-1}+\epsilon}.
\end{equation}
Then, by Lemma \ref{lem-montg}, (\ref{conj-formula}) and the above formulas,
\begin{eqnarray}
  &&\mathbf{P}\left(X(\boldsymbol{\theta})\geq c_1 V\right)\leq \mathbf{P}(X(\boldsymbol{\theta})\geq 2A(T))\leq \exp\left(-\frac{3}{4}A(T)^2B(T)^{-1}\right)\nonumber\\
&&\qquad\qquad\qquad~~\leq \exp\left(-c_2 V^2 T^{\frac{1}{k}-\epsilon}\right) \leq \exp\left(-c_3 V^{\frac{2k}{k-1}-\epsilon}\right).\nonumber
\end{eqnarray}

\subsection{The lower bound}
Consider the expectation of $e^{\lambda X(\boldsymbol{\theta})}$,
$\mathbf{E}(e^{\lambda X(\boldsymbol{\theta})})=\int e^{\lambda X(\boldsymbol{\theta})}d$.
By the definition of $X(\boldsymbol{\theta})$, we know that
\begin{equation}\nonumber
 \mathbf{E}(e^{\lambda X(\boldsymbol{\theta})})=\prod_{\gamma>0}I(\lambda r_{\gamma}),
\end{equation}
where $I(r)=\int_0^1 e^{r\sin 2\pi \theta}d\theta$.
Montgomery \cite{Mont} (formulas (2) and (3), pp.17) showed that
\begin{equation}\label{I-upper}
I(r)\leq \left\{
  \begin{array}{ll}
    e^r \\
    e^{\frac{r^2}{4}}
  \end{array}
\right. \quad\mbox{for all}\quad r\geq 0,
\end{equation}
and
\begin{equation}\label{I-lower}
  I(r)>\left\{
     \begin{array}{ll}
       2e^{\frac{r}{2}}, &r\geq 7, \\
       e^{\frac{r^2}{19}}, &0<r\leq 7.
     \end{array}
   \right.
\end{equation}

Now we want to find a $\lambda>0$ so that
\begin{equation}\label{E-equal}
  \mathbf{E}(e^{\lambda X(\boldsymbol{\theta})})=2\exp\left(\frac{\lambda}{2}\sum_{0<\gamma<T}r_{\gamma}\right).
\end{equation}
We show that such $\lambda$ exists. In fact, the two sides of (\ref{E-equal}) are continuous functions of $\lambda$; for $\lambda=0$ the left hand side is $1$ while the right hand side is $2$. Moreover, by (\ref{I-lower}), if $\lambda>\frac{7}{r_{\gamma}}$ for all $r_{\gamma}~(\gamma<T)$, then
\begin{equation}\nonumber
  \mathbf{E}(e^{\lambda X(\boldsymbol{\theta})})\geq \prod_{0<\gamma<T}I(\lambda r_{\gamma})\geq 2^{N(T)}\exp\left(\frac{\lambda}{2}\sum_{0<\gamma<T}r_{\gamma}\right)\geq 2\exp\left(\frac{\lambda}{2}\sum_{0<\gamma<T}r_{\gamma}\right),
\end{equation}
where $N(T)$ is the number of zeros of $\zeta(s)$ for $0<\gamma<T$. Thus, there is such a $\lambda$.
By (\ref{r-gamma-est}), we have, for $\gamma<T$,
$r_{\gamma}\gg \frac{1}{T^{\frac{1}{2}+\frac{1}{2k}+\epsilon}}$.
Thus,
\begin{equation}\label{lambda-bo}
  \lambda\leq c_4 T^{\frac{1}{2}+\frac{1}{2k}+\epsilon},
\end{equation}
for some constant $c_4>0$.

To get the lower bound, we need an inequality from \cite{Mont} (formula (4), pp.19), for any non-negative random variable $F$,
\begin{equation}\label{M-ineq}
  \mathbf{P}\left(F\geq \frac{1}{2}\mathbf{E}(F)\right)\geq \frac{\mathbf{E}(F)^2}{4\mathbf{E}(F^2)}.
\end{equation}
Let $F=e^{\lambda X(\boldsymbol{\theta})}$. By (\ref{I-upper}),
\begin{equation}\nonumber
 \mathbf{E}(F^2)\leq \exp\left(2\lambda\sum_{0<\gamma<T}r_{\gamma}+\lambda^2\sum_{\gamma\geq T}r^2_{\gamma}\right).
\end{equation}
So, by the above formula, (\ref{E-equal}), and (\ref{M-ineq}), we have
\begin{equation}\nonumber
  \mathbf{P}\left(X(\boldsymbol{\theta})\geq\frac{1}{2}\sum_{0<\gamma<T}r_{\gamma}\right)\geq \frac{\exp(\lambda\sum_{0<\gamma<T}r_{\gamma})}{\mathbf{E}(F^2)}\geq \exp\left(-\lambda\sum_{0<\gamma<T}r_{\gamma}-\lambda^2\sum_{\gamma\geq T}r^2_{\gamma}\right).
\end{equation}
Then, by (\ref{conj-formula}), (\ref{T-V}), and (\ref{lambda-bo}), we get
\begin{equation}\nonumber
  \mathbf{P}\left(X(\boldsymbol{\theta})\geq\frac{1}{2}\sum_{0<\gamma<T}r_{\gamma}\right)\geq \exp \left(-\lambda A(T)-\lambda^2 B(T)\right)\geq \exp\left(-c_5 T^{1+\epsilon}\right)\geq\exp\left(-c_6 V^{\frac{2k}{k-1}+\epsilon}\right).
\end{equation}

Hence, for any $\epsilon>0$,
\begin{equation}\nonumber
 \exp\left(-\widetilde{c}_1 V^{\frac{2k}{k-1}+\epsilon}\right)\leq \nu([V, \infty))\leq \exp\left(-\widetilde{c}_2 V^{\frac{2k}{k-1}-\epsilon}\right),
\end{equation}
for some constants $\widetilde{c}_1, \widetilde{c}_2 >0$ depending on $k$ and $\epsilon$.

\section{Large deviation conjecture and the order of $M_k(x)$}
In this section, we examine the tail of the limiting distribution more carefully, and heuristically derive a Large deviation conjecture and a conjecture about the maximum order of $M_k(x)$.

First, we give the proof of Theorem \ref{thm-la-dev-upper}. For any integer $k\geq 2$ and any integer $l\geq 1$, taking $w=\frac{1}{k}$ in Theorem \ref{thm-moments-zeta}, and by the functional equation, we get
\begin{equation}\label{aver-zeta-f}
  \sum_{0<\gamma\leq T}|\zeta(\frac{\rho}{k})|^{2l}\asymp T^{2l(\frac{1}{2}-\frac{1}{2k})}T\log T.
\end{equation}
By H\"{o}lder's inequality,
\begin{equation}\nonumber
  \sum_{0<\gamma\leq T}\frac{|\zeta(\frac{\rho}{k})|}{|\zeta'(\rho)|}\leq \left(\sum_{0<\gamma\leq T}|\zeta(\frac{\rho}{k})|^{2l}\right)^{\frac{1}{2l}}\left(\sum_{0<\gamma\leq T}\left|\frac{1}{\zeta'(\rho)}\right|^{\frac{2l}{2l-1}}\right)^{\frac{2l-1}{2l}}.
\end{equation}
Since $l\geq 1$, $\frac{l}{2l-1}<\frac{3}{2}$. Thus, by (\ref{aver-zeta-f}) and (\ref{G-conj}),
\begin{equation}\label{ladev-pf-up}
  \sum_{0<\gamma\leq T}\frac{|\zeta(\frac{\rho}{k})|}{|\zeta'(\rho)|}\ll_l T^{\frac{1}{2}-\frac{1}{2k}}\cdot T(\log T)^{\frac{l}{4l-1}}.
\end{equation}
Also, for any $0<\delta<1$, by H\"{o}lder's inequality, we have
\begin{equation}\label{lo-thm-devi}
\sum_{\delta T\leq\gamma\leq T}\frac{1}{|\zeta'(\rho)|^{\frac{2l}{2l+1}}}=\sum_{\delta T\leq\gamma\leq T}\left|\frac{\zeta(\frac{\rho}{k})}{\zeta'(\rho)}\right|^{\frac{2l}{2l+1}}\cdot \frac{1}{|\zeta(\frac{\rho}{k})|^{\frac{2l}{2l+1}}}\leq \left( \sum_{\delta T\leq\gamma\leq T}\frac{|\zeta(\frac{\rho}{k})|}{|\zeta'(\rho)|}\right)^{\frac{2l}{2l+1}}\left(\sum_{\delta T\leq\gamma\leq T}\frac{1}{|\zeta(\frac{\rho}{k})|^{2l}} \right)^{\frac{1}{2l+1}}.
\end{equation}
By (\ref{G-conj}), there exists a small enough $\delta>0$ such that,
$$\sum_{\delta T\leq\gamma\leq T}\frac{1}{|\zeta'(\rho)|^{\frac{2l}{2l+1}}}\asymp T(\log T)^{\left(\frac{l+1}{2l+1} \right)^2}.$$
For such $\delta$, by Theorem \ref{thm-moments-zeta} and the functional equation,
$$\sum_{\delta T\leq\gamma\leq T}\frac{1}{|\zeta(\frac{\rho}{k})|^{2l}}\asymp\frac{T\log T}{T^{(\frac{1}{2}-\frac{1}{2k})\cdot 2l}}.$$
So, by (\ref{lo-thm-devi}), we deduce that
\begin{equation}\label{ladevi-pf-lo}
 \sum_{0<\gamma\leq T}\frac{|\zeta(\frac{\rho}{k})|}{|\zeta'(\rho)|}\geq\sum_{\delta T\leq\gamma\leq T}\frac{|\zeta(\frac{\rho}{k})|}{|\zeta'(\rho)|}\gg_l T\cdot T^{\frac{1}{2}-\frac{1}{2k}}(\log T)^{\frac{l}{2(2l+1)}}.
\end{equation}
Thus, by (\ref{ladev-pf-up}) and (\ref{ladevi-pf-lo}), we get
\begin{equation}\nonumber
  T^{\frac{1}{2}-\frac{1}{2k}}\cdot T(\log T)^{\frac{1}{4}-o(1)}\ll\sum_{0<\gamma\leq T}\frac{|\zeta(\frac{\rho}{k})|}{|\zeta'(\rho)|}\ll T^{\frac{1}{2}-\frac{1}{2k}}\cdot T(\log T)^{\frac{1}{4}+o(1)}.
\end{equation}
Similarly, we have
\begin{equation}\nonumber
  \sum_{0<\gamma\leq T}\frac{|\zeta(\frac{\rho}{k})|^2}{|\zeta'(\rho)|^2}\ll_l T^{1-\frac{1}{k}}\cdot T (\log T)^{\frac{1}{l-1}}\ll T^{1-\frac{1}{k}}\cdot T (\log T)^{o(1)}.
\end{equation}
Then, by partial summation,
\begin{equation}\label{A-B-tr-result}
 T^{\frac{1}{2}-\frac{1}{2k}}\cdot (\log T)^{\frac{1}{4}-o(1)}\ll A(T)=\sum_{0<\gamma<T}\frac{2|\zeta(\frac{\rho}{k})|}{|\rho\zeta'(\rho)|}\ll T^{\frac{1}{2}-\frac{1}{2k}}(\log T)^{\frac{1}{4}+o(1)}, \
\end{equation}
and
\begin{equation}\label{AB-B-tr-result}
B(T)=\sum_{\gamma\geq T}\frac{4|\zeta(\frac{\rho}{k})|^2}{|\rho\zeta'(\rho)|^2}\ll \frac{(\log T)^{o(1)}}{T^{\frac{1}{k}}}.
\end{equation}

We choose $T$ such that $A(T^-)<V\leq A(T)$. Then, by (\ref{A-B-tr-result}),
\begin{equation}\nonumber
T^{\frac{1}{2}-\frac{1}{2k}}(\log T)^{\frac{1}{4}-o(1)}\ll V\ll T^{\frac{1}{2}-\frac{1}{2k}}(\log T)^{\frac{1}{4}+o(1)}.
\end{equation}
So, we have
\begin{equation}\nonumber
\left(\frac{V}{(\log V)^{\frac{1}{4}+o(1)}}\right)^{\frac{2k}{k-1}}\ll  T\ll \left(\frac{V}{(\log V)^{\frac{1}{4}-o(1)}}\right)^{\frac{2k}{k-1}}.
\end{equation}
Then, by (\ref{A-B-tr-result}), (\ref{AB-B-tr-result}), and Lemma \ref{lem-montg}, we  get (\ref{conj-up-bound}).
Hence, Theorem \ref{thm-la-dev-upper} follows.

By (\ref{A-B-tr-result}) and (\ref{AB-B-tr-result}), the conjectured formulas are
\begin{equation}\label{est-A-B}
A(T)\asymp T^{\frac{1}{2}-\frac{1}{2k}}(\log T)^{\frac{1}{4}}, \quad\mbox{and}\quad B(T)\asymp \frac{1}{T^{\frac{1}{k}}}.
\end{equation}
Then, by Lemma \ref{lem-montg}, we get
\begin{equation}\label{ladevi-up-bound}
  \nu([V, \infty))\ll \exp\left(-c_1' \frac{V^{\frac{2k}{k-1}}}{(\log V)^{\frac{1}{2(k-1)}}}\right).
\end{equation}
In Remark $2.4$ of \cite{HKC}, the authors mentioned that: if one redefines $J_{-l}(T)$ to exclude these rare points, where $|\zeta'(\frac{1}{2}+i\gamma_n)|$ is very close to zero, then the Random Matrix Theory should still predict the universal behavior. Thus, if we let $\{r_{\gamma'}\}$ be the decreasing sequence after reordering the sequence $\{r_{\gamma}\}$, we conjecture that we still have the similar estimates like (\ref{est-A-B}), i.e.
\begin{equation}\label{reorder-est}
 A'(T):=\sum_{0<\gamma'<T}r_{\gamma'}\asymp T^{\frac{1}{2}-\frac{1}{2k}}(\log T)^{\frac{1}{4}}, \quad\mbox{and}\quad B'(T):=\sum_{\gamma'\geq T}r^2_{\gamma'}\asymp \frac{1}{T^{\frac{1}{k}}}.
\end{equation}

Moreover, Hattori and Matsumoto (\cite{HM}, Theorem 4) proved an equivalent condition in terms of $A'$ and $B'$ for the existence of the lower bound of Montgomery type (the type of lower bound in Theorem 1 of \cite{Mont}). Hence, by Theorem 1 in \cite{Mont}, Theorem 4 in \cite{HM}, and (\ref{reorder-est}), we conjecture that the upper bound (\ref{ladevi-up-bound}) gives the correct order of $\nu([V, \infty))$, i.e. the Large deviation conjecture (\ref{conj-order}).

Then, using a similar heuristic analysis to section 4.3 of \cite{Ng}, we make the conjecture (\ref{M-conj}).

\section{Moments of $\zeta(1-w\rho)$ for $0<w<1$}\label{sec-zeta-moments}
In this section, we give the proof of Theorem \ref{thm-moments-zeta}.

Let $l\geq 1$ be an integer. Taking $c=\frac{1+w}{4w}>\frac{1}{2}$. By the residue theorem,
\begin{eqnarray}\label{thm-mom-pf-total}
 && \sum_{0<\gamma\leq T}|\zeta(1-w\rho)|^{2l}\nonumber\\
  &&\quad =\frac{1}{2\pi i}\left(\int_{c+i}^{c+iT}+\int_{c+iT}^{1-c+iT}+\int_{1-c+iT}^{1-c+i}+\int_{1-c+i}^{c+i}\right)\zeta^l(1-ws)\zeta^l(1-w(1-s))\frac{\zeta'}{\zeta}(s)ds\nonumber\\
  &&\quad =: J_1+J_2+J_3+J_4.
\end{eqnarray}
By \cite{Daven} (page 108), we may assume that $T$ satisfies
$$|\gamma-T|\gg \frac{1}{\log T}, \quad \mbox{for all ordinates}~\gamma$$
and
$$\frac{\zeta'}{\zeta}(\sigma+iT)\ll (\log T)^2, \quad \mbox{uniformly for all}~1-c\leq\sigma\leq c.$$
For general $T$, since $w<1$, under the Riemann Hypothesis, and by (\ref{zero-count}), the error is
$$\ll \sum_{T\leq \gamma\leq T+1}|\zeta(1-w\rho)|^{2l}\ll_{w, l} T^{\epsilon}.$$
Since $0<w<1$ and $c=\frac{1+w}{4w}$,  for $1-c\leq \sigma\leq c$,
$1-w\sigma\geq 1-wc=1-\frac{1+w}{4}>\frac{1}{2}\quad\mbox{and}\quad 1-w(1-\sigma)\geq 1-wc>\frac{1}{2}.$
So under the Riemann Hypothesis,
\begin{equation}\label{thm-mom-pf-J2}
  J_2=\frac{i}{2\pi}\int_{1-c}^{c}\zeta^l(1-w\sigma-iwT)\zeta^l(1-w(1-\sigma)+iwT)
  \frac{\zeta'}{\zeta}(\sigma+iT)d\sigma\ll_{w,l} T^{\epsilon}.
\end{equation}
Similarly, $J_4\ll_{w, l} 1$.
We relate $J_3$ to $J_1$,
\begin{eqnarray}
  J_3&=&\frac{1}{2\pi}\int_{T}^1 \zeta^l(1-w(1-c)-iwt)\zeta^l(1-wc+iwt)\frac{\zeta'}{\zeta}(1-c+it)dt\nonumber\\
  &=&-\frac{1}{2\pi}\overline{\int_1^T \zeta^l(1-w(1-c)+iwt)\zeta^l(1-wc-iwt)\frac{\zeta'}{\zeta}(1-c-it)dt}.\nonumber
\end{eqnarray}

By the functional equation $\zeta(s)=\chi(s)\zeta(1-s)$, we find that
\begin{eqnarray}
  J_3&=&-\frac{1}{2\pi}\overline{\int_1^T \zeta^l(1-w(1-c)+iwt)\zeta^l(1-wc-iwt)\frac{\chi'}{\chi}(1-c-it)dt}\nonumber\\
  &&+\frac{1}{2\pi}\overline{\int_1^T \zeta^l(1-w(1-c)+iwt)\zeta^l(1-wc-iwt)\frac{\zeta'}{\zeta}(c+it)dt},\nonumber
\end{eqnarray}
where
$\chi(s)=2^s\pi^{s-1}\Gamma(1-s)\sin(\frac{\pi s}{2}).$
By Stirling's formua,
\begin{equation}\nonumber
  -\frac{\chi'}{\chi}(1-c-it)=\log \left(\frac{|t|}{2\pi}\right)\left(1+O\left(\frac{1}{|t|}\right)\right).
\end{equation}
The term $O\left(\frac{1}{|t|}\right)$ contributes to $J_3$ an amount of $O_{w,l}(T)$. Let
\begin{equation}\nonumber
  K=\frac{1}{2\pi}\int_1^T \zeta^l(1-w(1-c)+iwt)\zeta^l(1-wc-iwt)\log(\frac{t}{2\pi})dt.
\end{equation}
Then,
\begin{equation}\label{thm-mom-pf-J3}
  J_3=K+\overline{J_1}+O_{w,l}(T).
\end{equation}
First, by the residue theorem, we calculate
\begin{eqnarray}
  I(T)&:=&\int_1^T \zeta^l(1-w(1-c)+iwt)\zeta^l(1-wc-iwt)dt\nonumber\\
  &=&\frac{1}{i}\int_{c+i}^{c+iT} \zeta^l(1-w(1-s))\zeta^l(1-ws)ds.\nonumber\\
&=&\frac{1}{i}\left(\int_{c+i}^{\frac{1}{2}+i}+\int_{\frac{1}{2}+i}^{\frac{1}{2}+iT}+\int_{\frac{1}{2}+iT}^{c+iT}\right)\zeta^l(1-w(1-s))\zeta^l(1-ws)ds\nonumber\\
  &=&\int_{1}^{T}\zeta^l(1-\frac{w}{2}+iwt)\zeta^l(1-\frac{w}{2}-iwt)dt+O_{w,l}(T^{\epsilon})\nonumber\\
  &=&\frac{1}{w}\int_{w}^{wT}\left|\zeta(1-\frac{w}{2}+it)\right|^{2l}dt+O_{w,l}(T^{\epsilon}).\nonumber
\end{eqnarray}
For $\frac{1}{2}<\sigma<1$ and $\frac{x}{2}<t<x$, taking $T=x^{3}$ in the proof of Theorem 13.3 in \cite{Titch} (page 330), we get
\begin{equation}\nonumber
  \zeta^l(s)=\sum_{n<x}\frac{d_l(n)}{n^s}+O(x^{-\epsilon}).
\end{equation}
Then, by Montgomery and Vaughan's mean value theorem for Dirichlet polynomials (Lemma 1 in \cite{Tsang}, originally due to \cite{MontVaug}, Corollary 3 to Theorem 2), we deduce that
\begin{equation}\label{appe-mean-zeta}
  I(T)= \left(\sum_{n=1}^{\infty}\frac{d_l^2(n)}{n^{2-w}}\right) T+O_{w,l}(T^{1-\epsilon}).
\end{equation}
By partial summation, we get
\begin{equation}\label{thm-mom-pf-K}
  K=\frac{1}{2\pi}\log(\frac{T}{2\pi})I(T)-\frac{1}{2\pi}\int_1^T\frac{I(t)}{t}dt=\frac{1}{2\pi}\left(\sum_{n=1}^{\infty}\frac{d_l^2(n)}{n^{2-w}}\right) T\log T+O_{w,l}(T).
\end{equation}

Now, we estimate $J_1$. By Theorem 14.5 in \cite{Titch} (page 341),
\begin{equation}\nonumber
  \frac{\zeta'}{\zeta}(s)=O((\log t)^{2-2\sigma})\quad \mbox{uniformly for}\quad \frac{1}{2}<\sigma_0\leq \sigma\leq\sigma_1<1.
\end{equation}
Thus, by (\ref{appe-mean-zeta}) and H\"{o}lder's inequality,
\begin{equation}\label{thm-mom-pf-J1}
  |J_1|\ll_{w,l} T(\log T)^{\frac{1}{2}}.
\end{equation}

Combining (\ref{thm-mom-pf-total}), (\ref{thm-mom-pf-J2}), (\ref{thm-mom-pf-J3}), (\ref{thm-mom-pf-K}), and (\ref{thm-mom-pf-J1}), we get (\ref{thm-moms-one}). By Theorem 14.25 (A) in \cite{Titch}, under the Riemann Hypothesis, the series $\sum_{n=1}^{\infty}\frac{\mu(n)}{n^s}$ is convergent, and its sum is $\frac{1}{\zeta(s)}$, for every $s$ with $\sigma>\frac{1}{2}$. Hence, the same proof also works for negative powers, and we have (\ref{thm-moms-two}). \qed

\vspace{1em}
\textbf{Acknowledgements.} This paper was finished under the support of the NSF grant DMS-1201442. I would like to thank my advisor, Prof. Kevin Ford, for providing me the financial support, giving me invaluable suggestions and encouragements to finish this project and very useful comments. The author is also grateful to the referee for suggesting the reference of $B^2$-almost periodic functions and helpful comments.


\vspace{2em}

DEPARTMENT OF MATHEMATICS, UNIVERSITY OF ILLINOIS AT URBANA-CHAMPAIGN, 1409 WEST GREEN ST., URBANA, IL 61801, USA\\
\emph{E-mail address}: xmeng13@illinois.edu

\end{document}